\documentclass[12pt]{article}
\usepackage[final]{epsfig}
\usepackage{graphics}
\usepackage{amsmath}
\usepackage{amsfonts}
\usepackage{latexsym}
\usepackage{amssymb}
\usepackage{graphicx}
\usepackage{epstopdf}
\usepackage{url}

\newtheorem{lemma}{Lemma}[section]
\newtheorem{proposition}[lemma]{Proposition}
\newtheorem{remark}[lemma]{Remark}
\newtheorem{example}[lemma]{Example}
\newtheorem{theorem}{Theorem}
\newtheorem{definition}{Definition}
\newtheorem{corollary}[lemma]{Corollary}

\newtheorem{conjecture}{Conjecture}

\begin{document}
\newcommand{\eps}{{\varepsilon}}
\newcommand{\g}{{\gamma}}
\newcommand{\G}{{\Gamma}}
\newcommand{\proofend}{$\Box$\bigskip}
\newcommand{\C}{{\mathbf C}}
\newcommand{\Q}{{\mathbf Q}}
\newcommand{\R}{{\mathbf R}}
\newcommand{\Z}{{\mathbf Z}}
\newcommand{\RP}{{\mathbf {RP}}}
\newcommand{\CP}{{\mathbf {CP}}}
\newcommand{\Tr}{{\rm Tr\ }}
\def\proof{\paragraph{Proof.}}

\title {Configuration spaces of plane polygons and a sub-Riemannian approach to the equitangent problem}
\author{Jes\'us Jer\'onimo-Castro\footnote{
Facultad de Ingenier\'ia, Universidad Aut\'onoma de Quer\'etaro, Cerro de las Campanas s/n, C. P. 76010, Quer\'etaro, M\'exico; jesusjero@hotmail.com
}
\ and Serge Tabachnikov\footnote{
Department of Mathematics,
Pennsylvania State University,
University Park, PA 16802, USA,
and ICERM, ICERM, Brown University, Box 1995, Providence, RI 02912, USA;
tabachni@math.psu.edu} }
\maketitle

\section{Introduction} \label{intro}

\begin{figure}[hbtp] 
\centering
\includegraphics[height=1.7in]{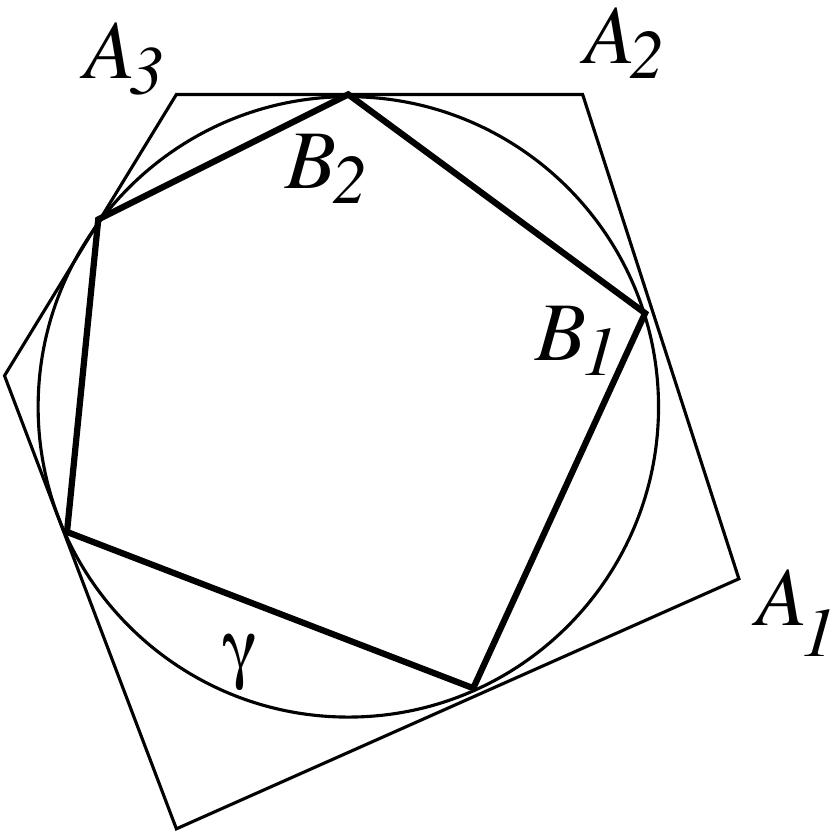}
\caption{An equitangent pentagon: $|A_2B_1|=|A_2 B_2|$, and likewise, cyclically.}
\label{framed1}
\end{figure}

Given a strictly convex plane curve $\gamma$ and a point $x$ in
its exterior, there are two tangent segments to $\gamma$ from $x$.
The {\it equitangent locus} of $\gamma$ is the set of points $x$
for which the two tangent segments have equal lengths. An {\it equitangent $n$-gon} 
of $\gamma$ is a circumscribed $n$-gon such that the tangent segments to $\gamma$ from the vertices have equal lengths, see Figure \ref{framed1}.

For example, if $\gamma$ is a circle, the  equitangent locus is
the whole exterior of $\gamma$ (this property is characteristic of
circles \cite{RT}), and if $\gamma$ is an ellipse then its
equitangent locus consists of the two axes of symmetry (this implies that an ellipse does not admit equitangent $n$-gons with $n\neq 4$).

Generically, the equitangent locus of a curve $\gamma$ is also a
curve, say, $\Gamma$.  Note that $\Gamma$ is not empty; in fact,
every curve $\gamma$ has pairs of equal tangent segments  of an
arbitrary length. The equitangent problem is to study the relation
between a curve and its equitangent locus.

For example, if the equitangent locus  contains a line tangent to $\gamma$, then
$\gamma$ must be a circle. On the other hand, there exists an infinite-dimensional,  functional
family of curves $\g$ for which the equitangent locus contains a line $\G$, disjoint from $\g$: these 
are the  curves of constant width in the
half-plane model of the hyperbolic geometry, where $\Gamma$ plays
the role of the absolute. These observations were made in
\cite{JR,JRT}.

Another, somewhat surprising, observation is that there exist
nested curves $\gamma \subset \gamma_1$ such that  $\gamma_1$ is
disjoint from the equitangent locus of $\gamma$, see \cite{Ta2}.
In other words, a point can make a full circuit around $\gamma$ in
such a way that, at all times, the two tangent segments to $\g$ have
unequal lengths.

\begin{figure}[hbtp]
\centering
\includegraphics[height=1.2in]{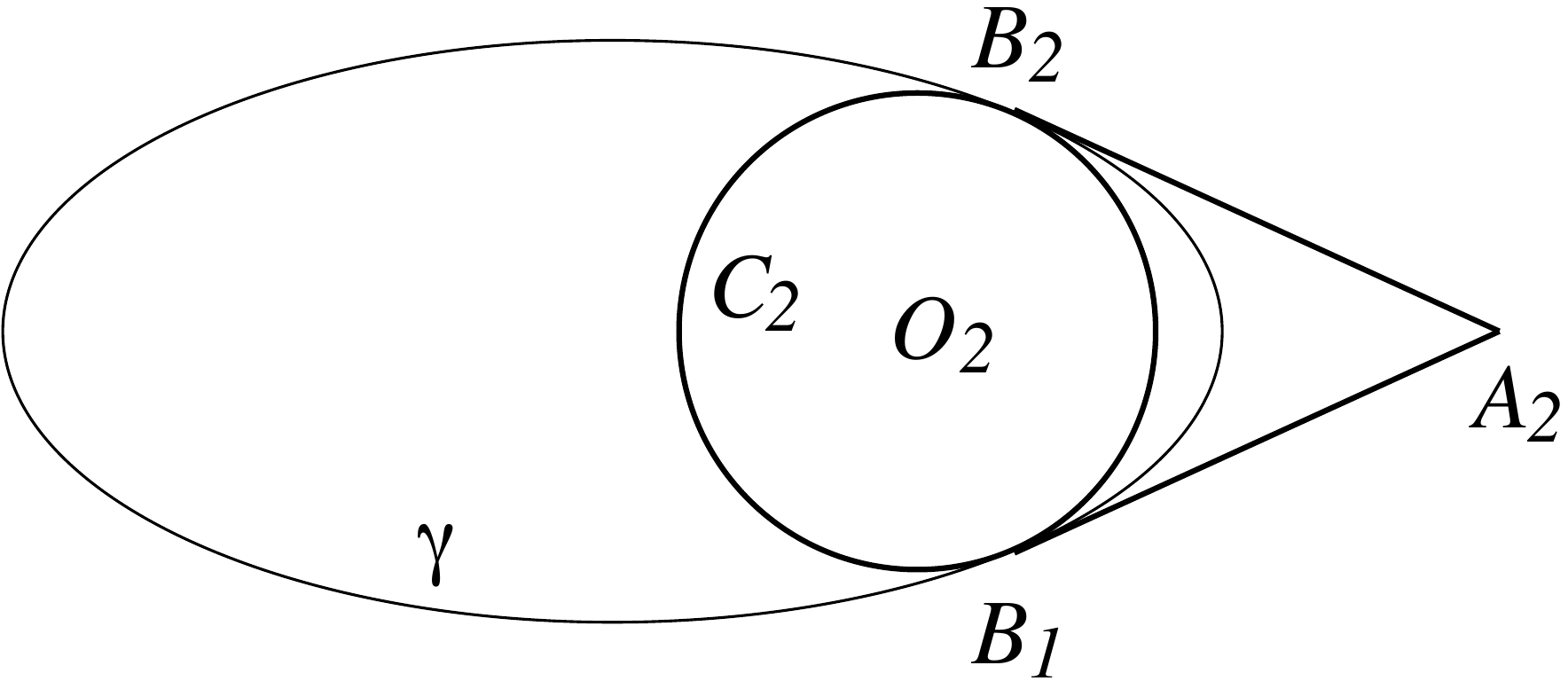}
\caption{$C_2$ is a bitangent circle, its center $O_2$ belongs to the symmetry set of the 
curve $\gamma$.}
\label{ellipse}
\end{figure}

Two tangent segments to a curve $\gamma$ are equal if and only if
there is a circle touching $\gamma$ at these tangency points, see Figure \ref{ellipse}. The
locus of centers of such bitangent circles is called the {\it
symmetry set} of $\gamma$, see \cite{BGG, GB}. Symmetry sets are
of great interest,  due to their applications in image recognition
and computer vision.

Our approach to the equitangent problem uses ideas of sub-Riemannian geometry.

Yu. Baryshnikov and V. Zharnitsky \cite{BZ}, and J. Landsberg
\cite{La}, pioneered a sub-Riemannian approach in the study of
mathematical billiards. One is interested in plane billiard tables
that admit a 1-parameter family of $n$-periodic billiard
trajectories. If the billiard curve is an ellipse, such a family
exists for every  $n\ge 3$. Are there other curves with this
property? (For $n=2$, the billiard table has constant width).

It is observed in \cite{BZ,La} that if an $n$-gon is a periodic
billiard trajectory then the directions of the boundary of the
billiard table at its vertices are determined by the reflection law (the angle of incidence equals the angle of reflection). This defines an
$n$-dimensional distribution, called the {\it Birkhoff
distribution}, on the space of $n$-gons with a fixed perimeter length, and a 1-parameter family
of  periodic billiard trajectories corresponds to a curve in this
space, tangent to the distribution (a horizontal curve).

This point of view makes it possible to construct a  functional
space of billiard tables, admitting 1-parameter families of
$n$-periodic billiard trajectories, and to give a new proof of the $n=3$ case of the Ivrii conjecture (the set of $n$-periodic billiard trajectories has zero measure). 
See  \cite{GT,TZ} for a similar approach to outer billiards.

Let us outline the contents of the present paper. 
We are interested in the situation when a convex curve $\gamma$ admits a 1-parameter family of equitangent $n$-gons; such curves are called {\it $n$-fine curves}, see Section \ref{distrF} for the precise definition. 

Consider Figure \ref{framed1} again. The  vectors $B_1 A_2$ and $B_2 A_3$ make equal angles with the vector $B_1 B_2$. A polygon $B_1 B_2 \ldots$ is called {\it framed} if vectors are assigned to its vertices such that the vectors at the adjacent vertices make equal angles with the respective side; see Section \ref{chains} for  definition. The space of framed $n$-gons is denoted by ${\cal F}_n$.

We prove that ${\cal F}_n$ is an $2n$-dimensional variety. For $n$ odd, an $n$-gon uniquely determines the framing (up to a certain involution, see Section \ref{chains}). If $n\ge 4$ is even, a framed polygon satisfies one non-trivial condition, and if this condition  holds then there exists a 1-parameter family of framings. For example, a framed quadrilateral is cyclic (inscribed into a circle).

If $\gamma$ is an $n$-fine curve then one can turn the respective framed $n$-gon $B_1 B_2 \ldots$ inside $\gamma$; in this motion, the velocity of each vertex is aligned with the respective framing vector. This condition defines an $n$-dimensional distribution on the space ${\cal F}_n$, denoted by ${\cal D}$. The 1-parameter family of the framed polygons $B_1 B_2 \ldots$ corresponds to a horizontal curve in $({\cal F}_n, {\cal D})$.

We prove that, on the open dense subset of {\it generic} framed polygons (see Section \ref{chains} for definition), the distribution ${\cal D}$ is bracket generating of the type $(n,2n)$, that is, the vector fields tangent to the
distribution and their first  commutators  span the tangent space to the space of framed polygons at every point.

Return to Figure \ref{ellipse}. If $A_1 A_2 \ldots$ is an equitangent $n$-gon then one obtains a {\ chain of circles} $C_1 C_2 \ldots$ in which every two consecutive circles are tangent. This chain is {\it oriented} in the sense that the circles  can rotate, as if they were linked gears. Let ${\cal C}_n$ be the space of such oriented chains of $n$ circles. This space also carries an $n$-dimensional distribution, ${\cal \widetilde B}$, see Section \ref{chains}, and an $n$-fine curve corresponds to a horizontal curve in $({\cal C}_n, {\cal \widetilde B})$.

We define an open dense subset of generic chains of circles, and prove that this set is in bijection with the set of  generic framed $n$-gons; this bijection identifies the distributions ${\cal D}$ and ${\cal \widetilde B}$.

The centers $O_1 O_2 \ldots$ of the circles in an oriented chain form an $n$-gon. The sides of this polygon are oriented, see Section \ref{chains}, and its signed perimeter length is equal to zero. This polygon is an {\it evolute}  of the respective piece-wise circular curve, see \cite{BG}, and the fact that the signed perimeter length vanishes is analogous to the well-known property of the evolute of a smooth closed curve: its signed length equals zero (the sign changes at the cusps). The space of  such zero-length $n$-gons is denoted by ${\cal E}_n$. The variety ${\cal E}_n$ has dimension $2n-1$, and  the projection ${\cal C}_n \to {\cal E}_n$ is an $\R$-bundle.

The projection of the distribution ${\cal \widetilde B}$ is an $n$-dimensional distribution ${\cal B}$ on ${\cal E}_n$, essentially, the above mentioned Birkhoff distribution, defined by the billiard reflection law. This provides a somewhat  unexpected connection between the equitangent and the billiard problems. 

The case $n=2$, that of `bigons', is special, and we study it in Section \ref{chord}. The space ${\cal F}_2$ is 5-dimensional, and the distribution ${\cal D}$ is bracket generating of the type $(3,5)$. We show that, given a strictly convex $2$-fine curve $\g$, there exists a set of $2$-fine curves sufficiently close to $\g$, parameterized by a certain infinite-dimensional Hilbert manifold. Thus $2$-fine curves enjoy a great flexibility. 

Section \ref{tri} concerns the cases $n=3$ and $n=4$: we prove that, in these cases, one has rigidity:  $3$- and $4$-fine curves must be circles. 

Section \ref{exa} provides a construction of nested curves $\gamma \subset \Gamma$ where $\G$ is contained in the equitangent locus of $\g$. In other words, the two tangent segments to $\g$ from every point of $\G$ have equal lengths. The curve $\g$ is a smooth strictly convex perturbation of a regular polygon.

In the last Section \ref{canit}, we address the following question: does there exist a non-circular $n$-fine curve such that the vertices of the respective equitangent $n$-gons belong to a circle? 

By the Poncelet Porism (see, e.g., \cite{Fl}), there exist  1-parameter families of bicentric polygons, inscribed into a circle $\G$ and circumscribed about a  circle $\g$, and these polygons are automatically equitangent to $\g$; we ask whether $\g$ can be a non-circular oval. We conjecture that the answer is negative, and establish an infinitesimal version of this conjecture for odd $n$. For odd $n$, we reduce this problem to the study of a certain flow on the space of inscribed $n$-gons.

\bigskip

{\bf Acknowledgements}. We are grateful to R. Montgomery and V. Zharnitsky for stimulating discussions. 
The second author was supported by the NSF grant DMS-1105442.


\section{Framed polygons and oriented chains of circles} \label{chains}

By the angle between two vectors $u$ and $v$ we mean the angle
through which one needs to rotate $u$ to align it with $v$ (with
the usual counter-clock orientation of the plane). We denote the
angle by $\angle (u,v)$. By this definition, $\angle (u,v) = -
\angle (v,u)$. The angles are measured mod $2\pi$.

The next definition is motivated by a discussion in Introduction, see Figure \ref{framed1}. 

\begin{definition} \label{framedpol}
A framing of an $n$-gon $B_1,B_2,\ldots$ is an assignment of unit
vectors $u_1,u_2,\ldots$ to its respective vertices such that
$$
\angle (u_i, B_i B_{i+1}) = \angle (B_i B_{i+1}, u_{i+1})
$$
for all $i=1,\ldots,n$ (here and elsewhere the index is understood
cyclically). If $(u_i)$ is a framing then so is $(-u_i)$, and we
consider framing up to this involution, acting simultaneously on
all vectors $u_i$.
\end{definition}

Denote by ${\cal F}_n$ the space of framed $n$-gons (we assume
that the polygons are non-degenerate in the sense that no two
consecutive vertices coincide).

Given a polygon $B$, denote by $\varphi_i$ the direction of the
side $B_i B_{i+1}$, and let $\theta_i=\varphi_i - \varphi_{i-1}$
be the exterior angle at vertex $B_i$.

\begin{lemma} \label{oddeven}
If $n$ is odd, then an $n$-gon has a unique framing. If $n$ is
even, then an $n$-gon admits a framing if and only if
$$
\sum_{i\ {\rm odd}} \theta_i \equiv \sum _{i\ {\rm even}} \theta_i \equiv 0\ {\rm mod}\ \pi,
$$
in which case there is a 1-parameter family of framings.
\end{lemma}

\proof Denote by $\alpha_i$ the direction of the framing vector
$u_i$. Then one has $ \alpha_i-\varphi_i=\varphi_i-\alpha_{i+1},\
i=1,\ldots,n. $ If $n$ is odd, the system has a unique solution
(up to addition of $\pi$ to all $\alpha_i$):
$$
\alpha_i = \sum_{j=0}^{n-1} (-1)^j \varphi_{i+j}.
$$

If $n$ is even then the system has corank 1, and taking the
alternating sum of the equations yields $2 \sum (-1)^i \varphi_i
=0$ mod $2\pi$. It follows that
$$
\sum_{i\ {\rm odd}} \theta_i\ \ {\rm and}\ \ \sum _{i\ {\rm even}} \theta_i
$$
are multiples of $\pi$.

If this condition holds then the space of solutions is
1-dimensional, and all solutions are obtained from one by adding a
constant to even angles $\alpha_i$ and subtracting this constant
from the odd ones. \proofend

For example, a framed quadrilateral is cyclic, that is, inscribed into a circle.

\begin{corollary} \label{dimF}
For $n\geq 3$, one has: dim ${\cal F}_n = 2n$.
\end{corollary}

Let $A_1, A_2, \ldots$ be an equitangent polygon circumscribed about
a curve $\gamma$, and let $B_1, B_2, \dots$ be the tangency points. Then, for
each $i$, there exists a circle, $C_i$, tangent to $\gamma$ at
points $B_{i-1}$ and $B_i$, see Figure \ref{ellipse}.

We obtain a chain of circles in which each next one is tangent to
the previous one (at point $B_i$). This prompts the next
definition.

\begin {definition} \label{chain}
A chain of circles $C_1,C_2,\ldots,C_n$ is a collection of circles
such that $C_i$ is tangent to $C_{i+1}$ for $i=1,\ldots, n$ (we
assume that $C_i \neq C_{i+1}$). A chain is called oriented if the
circles can be  assigned signs, subject to the following rule: the
signs of exterior tangent circles  are opposite, and the signs of
interior tangent circles  are the same. The signs in an oriented
chain can be simultaneously changed to the opposite; we consider
oriented chains up to this involution. A signed radius of a circle
is its radius, taken with the respective sign.
\end{definition}

The signs determine the coorientations of the circles; in an
oriented chain, the coorientations agree. In other words, the
circles in an oriented chain can rotate (as if they were linked
gears): circles of the same sign revolve in the same sense.

Note also that if a chain of $n$ circles is not oriented, one can
double it and obtain an oriented chain of $2n$ circles.

Denote by ${\cal C}_n$ the space of oriented chains of $n$ circles.

The next lemma is a version of Huygens' principle.

\begin{figure}[hbtp]
\centering
\includegraphics[height=1.7in]{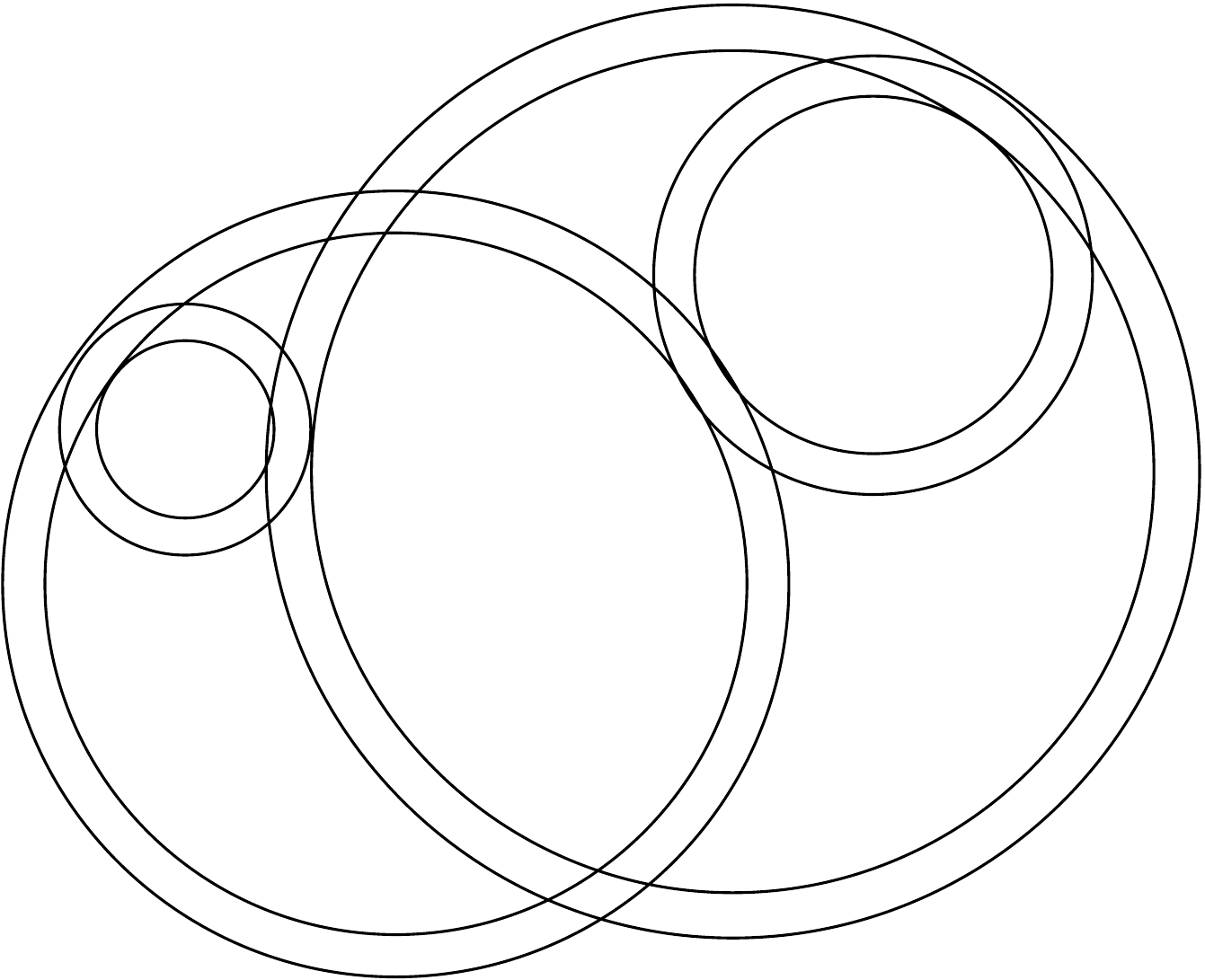}
\caption{An oriented chain of circles.}
\label{rot}
\end{figure}

\begin{lemma} \label{radius}
Given an oriented chain of circles, one can add a constant to all
the signed radii to obtain a new oriented chain, see Figure
\ref{rot}.
\end{lemma}

\proof This follows from the fact that the coorientations of the
circles in an oriented chain agree. \proofend

Let $O_i$ be the center of circle $C_i$ in an oriented chain.
Orient the segment $O_i O_{i+1}$ from the smaller (signed) radius
to the larger one (note that $r_i \neq r_{i+1}$, unless the two
circles coincide). The signed perimeter length of the polygon
$O_1, O_2,\ldots$ is the algebraic sum of the length of its sides,
where the side is taken with the positive sign if its orientation
agrees with the cyclic orientation of the polygon, and with the
negative sign otherwise.

\begin{lemma} \label{signed}
1) If a segment $O_1 O_2$ is oriented from $O_1$ to $O_2$ then its length is $r_2-r_1$, where the radii are signed.\\
2) The signed perimeter length of the polygon $O_1, O_2,\ldots, O_n$ is zero.
\end{lemma}

\begin{figure}[hbtp]
\centering
\includegraphics[height=0.7in]{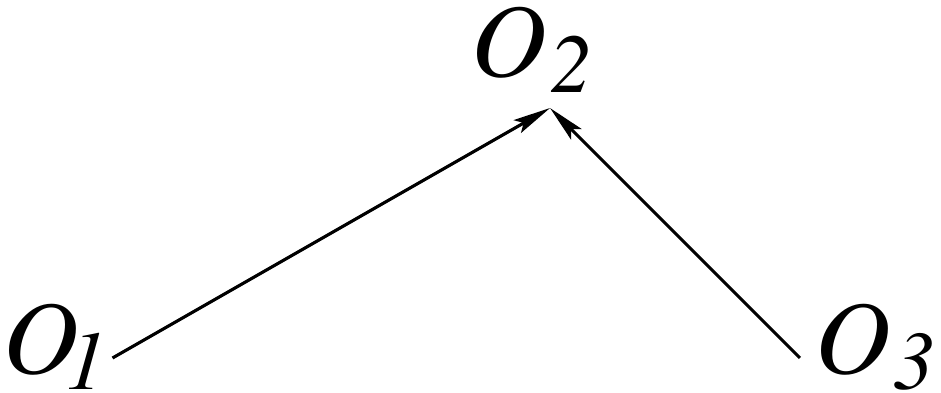}
\caption{Cancellation in the signed perimeter length: $(r_2-r_1)-(r_2-r_3)=r_3-r_1$.}
\label{tele}
\end{figure}

\proof Claim 1) easily follows by considering the three cases of
signs: $(-,-), (-,+), (+,+)$. For 2), taking the signs into
account, the sum `telescopes' to zero, see Figure \ref{tele}.
\proofend

Denote by ${\cal E}_n$ the space of $n$-gons with oriented edges
and zero signed perimeter length. One has the projection $\pi:
{\cal C}_n \to {\cal E}_n$ that assigns to a chain of circles the
polygon made of their centers.

\begin{lemma} \label{recons}
The projection $\pi: {\cal C}_n \to {\cal E}_n$ is an $\R$-bundle with the $\R$-action described in Lemma \ref{radius}.
\end{lemma}

\proof Consider a polygon with zero signed perimeter length, and
choose a vertex. Assign a signed radius to this vertex
arbitrarily, and consider an adjacent vertex. Using the first
claim of Lemma \ref{signed}, we assign the signed radius to this
vertex, and then proceed to the next one. Continuing in this way,
signed radii are assigned to all the vertices, and after we return
to the original vertex, this assignment is consistent, due to zero
signed perimeter length of the polygon. \proofend

The previous lemmas imply the next corollary.

\begin{corollary} \label{dimC}
One has: dim ${\cal C}_n=2n$.
\end{corollary}

Let $(B_i,u_i)$ be a framed polygon. It is called {\it
non-generic} if, for some $i$, the vectors $u_{i-1}, u_i, u_{i+1}$
are tangent to the circle through the vertices $B_{i-1}, B_i,
B_{i+1}$. A chain of circles $C_i$ is called {\it non-generic} if
the tangency points of $C_{i-1}$ with $C_i$ and of $C_i$ with
$C_{i+1}$ coincide.  Denote the spaces of generic chains and
generic framed polygons by ${\cal C}_n^\circ$ and ${\cal
F}_n^\circ$, respectively.

\begin{figure}[hbtp]
\centering
\includegraphics[height=2in]{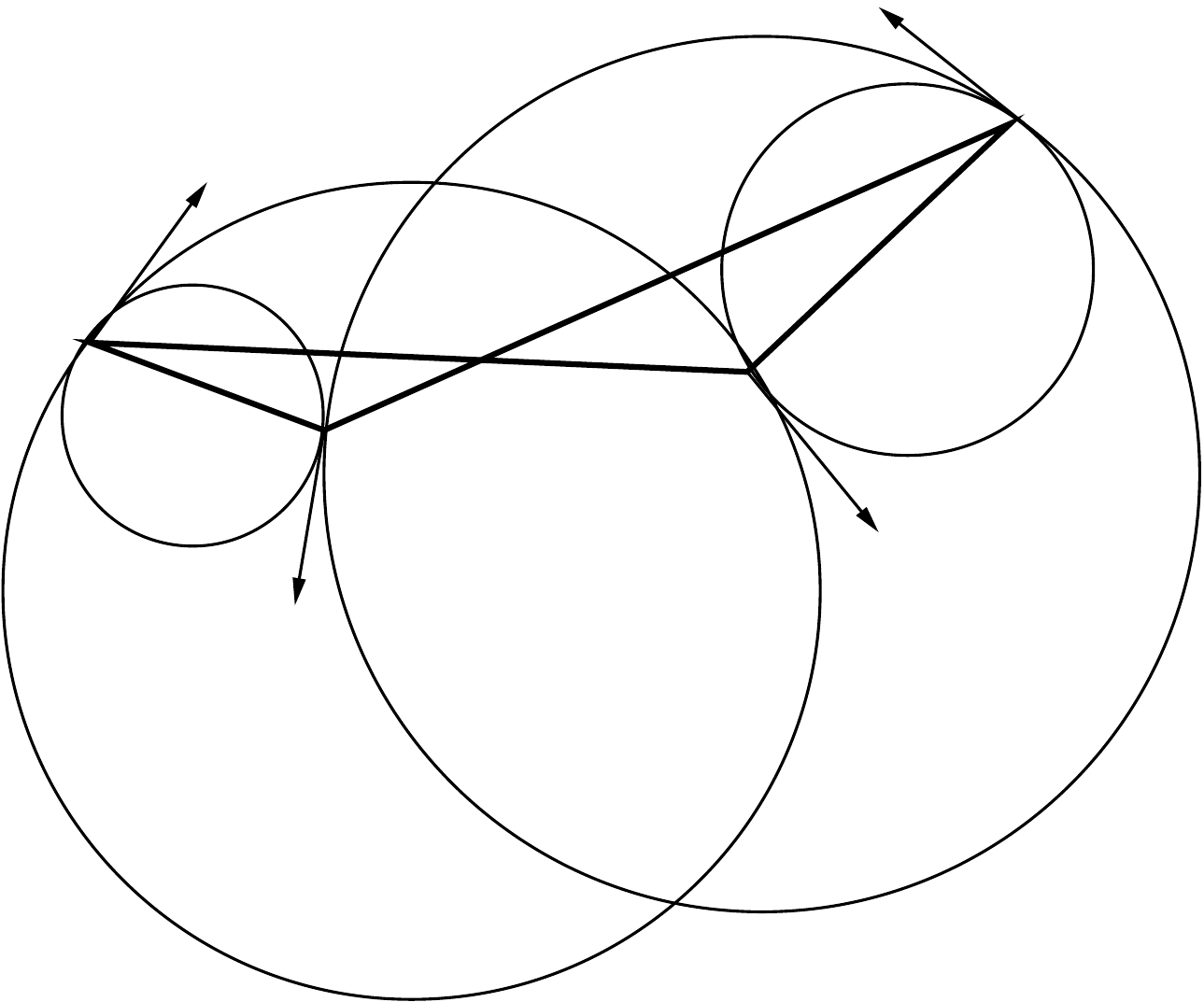}
\caption{The bijection ${\cal F}_n^\circ \to {\cal C}_n^\circ$.}
\label{bijec}
\end{figure}

The equality of the dimensions of ${\cal F}_n$ and ${\cal C}_n$
is not a coincidence. Let $B_{i-1}$ and $B_i$ be consecutive vertices of a framed polygon. There exists a unique circle through points $B_{i-1}$ and $B_i$, tangent to the framing vectors therein. This defines a map ${\cal F}_n^\circ \to {\cal C}_n^\circ$, see Figure \ref{bijec}.

\begin{proposition} \label{isom}
This map ${\cal F}_n^\circ \to {\cal C}_n^\circ$ is a bijection.
\end{proposition}

\proof A framed polygons gives rise to a chain of circles whose
consecutive tangency points are the vertices of the polygon. The
framing vectors define the directions of consistent rotation of
these circles, hence the chain is oriented. Since the framed polygon is generic, no two consecutive circles coincide.

Conversely, the  tangency points of the consecutive circles in a
chain define a polygon, and a choice of consistent rotations of
the circles provides its framing (by the velocity vectors of the
rotating circles). \proofend


\section{A distribution on ${\cal C}_n$} \label{distrC}

The Birkhoff distribution on the space of
polygons, introduced in \cite{BZ,La}, is defined as follows.

\begin{figure}[hbtp]
\centering
\includegraphics[height=0.9in]{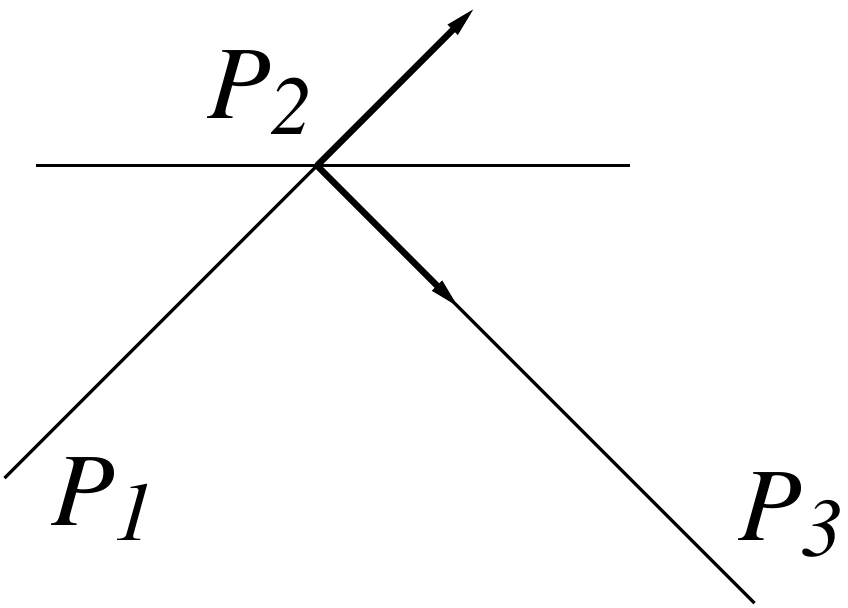}
\caption{Birkhoff distribution.}
\label{mirr}
\end{figure}

Let $P_1,P_2,\ldots$ be a polygon. Assign the direction at vertex
$P_i$ generated by the vector
$$
\frac{P_i - P_{i-1}}{|P_i - P_{i-1}|} + \frac{P_{i+1}-P_i}{|P_{i+1}-P_i|},
$$
see Figure \ref{mirr}. This agrees with the billiard reflection law:  the incoming trajectory $P_{i-1}
P_i$ elastically reflects to the outgoing trajectory $P_i P_{i+1}$.

Assigning a direction to every vertex defines an $n$-dimensional
distribution ${\cal B}$ on the $2n$-dimensional space of $n$-gons,
called the Birkhoff distribution. In fact, ${\cal B}$ is tangent
to the level hypersurfaces of the perimeter length function,
defining a distribution on these hypersurfaces.

It is proved in \cite{BZ} that ${\cal B}$ is bracket generating of
the type $(n, 2n -1)$, i.e., the vector fields tangent to the
distribution and their first  commutators already span the tangent
space to the space of $n$-gons of a fixed perimeter length.

We shall now construct a distribution on the space of oriented
chains of circles ${\cal C}_n$ which is an extension of the
Birkhoff distribution.

\begin{figure}[hbtp]
\centering
\includegraphics[height=0.9in]{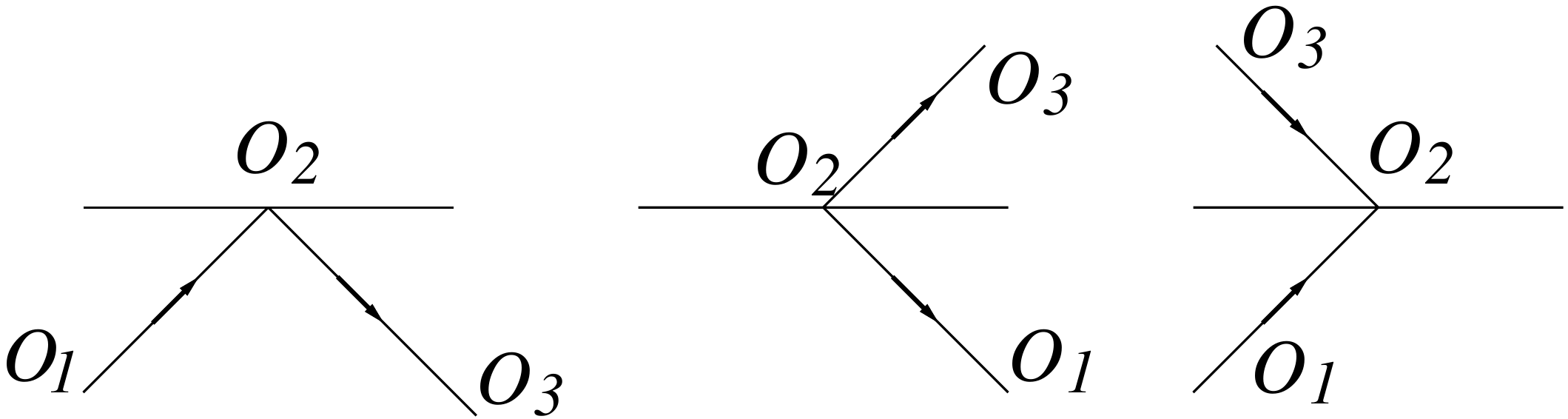}
\caption{Generalized billiard reflections.}
\label{3refl}
\end{figure}

Consider an oriented chain of circles, and let $O_1,O_2,\ldots$
and $r_1,r_2,\ldots$ be the centers and the signed radii. Let us
define an infinitesimal deformation of vertex $O_i$ and radius
$r_i$. The two sides, adjacent to $O_i$, are oriented. Consider
the line that bisects the angle between these two oriented lines:
this bisector defines infinitesimal motions of the vertex  $O_i$,
see Figure \ref{3refl}.\footnote{This kind of generalized billiards were recently studied by A. Glutsyuk \cite{Gl1,Gl2}.}
 The radius $r_i$ infinitesimally changes
in such a way that claim 1) of Lemma \ref{signed} continues to
hold.

The next lemma shows that the definition is correct. Let $O_i'$ be
the perturbed vertex.

\begin{lemma} \label{corr}
In the linear approximation, the signed perimeter length of the
polygon $O_1,\ldots, O_{i-1}, O_i', O_{i+1},\ldots, O_n$ remains
zero.
\end{lemma}

\proof This follows from the optical properties of conics.

Specifically, consider the first case in Figure \ref{3refl}. The
locus of points $O_2$ such that $|O_1 O_2|+|O_2 O_3|$ is fixed is
an ellipse with  foci $O_1$ and $O_3$. A billiard trajectory from
one focus reflects to another focus, i.e., the line in the figure
is tangent to the ellipse. Hence, in the first approximation,
$|O_1 O_2'|+|O_2' O_3|$ does not change.

The other two cases are similar: the locus of points $O_2$ such
that $|O_1 O_2|-|O_2 O_3|$ is fixed is a hyperbola with the foci
$O_1$ and $O_3$, and the claim follows from the optical property
of  hyperbolas. \proofend

We have defined an infinitesimal deformation for each $i$; taken
together, this gives an $n$-dimensional distribution on ${\cal
C}_n$ which we denote by ${\cal \widetilde B}$. Note that the rule
that assigns directions to the vertices of a polygon in ${\cal
E}_n$ defines  an $n$-dimensional distribution on ${\cal E}_n$.
Slightly abusing notation, we continue to denote it by ${\cal B}$
and call the Birkhoff distribution. The differential $d\pi$ is a
linear isomorphism of the spaces of the distributions ${\cal
\widetilde B}$ and ${\cal B}$.

The main result of this section is the following theorem.

\begin{theorem} \label{gen}
The distribution ${\cal \widetilde B}$ on ${\cal C}_n$ is bracket
generating of the type $(n,2n)$.
\end{theorem}

\proof
The proof is a modification of the argument given in \cite{BZ}.

Consider an oriented chain of circles, and let $O_1,O_2,\ldots$
and $r_1,r_2,\ldots$ be the centers and the signed radii. Let us
label the edges of the polygon $O_1,O_2,\ldots$ by half-integers
(also cyclic mod $n$). As before, the edges are oriented. Let
$u_j$ be the unit vector giving $j$th edge its orientation (so $j
\in \Z_n+1/2$), and let $2 \theta_i = \angle (u_{i-1/2},
u_{i+1/2}).$

Denote the rotation of the plane through $\pi/2$ by $J$ (the
complex structure on $\R^2=\C$). A tangent vector to ${\cal C}_n$
defines an infinitesimal motion of the vertices $O_i$ and a rate
of change of the radii $r_i$. Let $v_i$ be the  tangent vector to
${\cal C}_n$ that moves  vertex $O_i$ only with the velocity
$J(u_{i+1/2} - u_{i-1/2})$ (all other vertices do not move) and
changes  the radius $r_i$ only with the rate $-\sin(2\theta_i)$
(all other radii are intact). Obviously, the vectors $v_i$ are
linearly independent.

\begin{lemma} \label{span}
The vectors $v_i,\ i=1,\ldots,n$, span ${\cal \widetilde B}$.
\end{lemma}

\proof It is clear that the vector $J(u_{i+1/2} - u_{i-1/2})$
defines the `correct' direction at $O_i$, see Figure \ref{3refl}.
What one needs to check is that $r_i$ changes as claimed. To this
end, we use Lemma \ref{signed}.

Assume that the edge $O_{i-1} O_i$ is oriented toward $O_i$ (the
other case is analogous). Then $|O_{i-1} O_i|=r_i-r_{i-1}$. The
rate of change of $|O_{i-1} O_i|$ is
$$
u_{i-1/2} \cdot J(u_{i+1/2} - u_{i-1/2}) = (u_{i+1/2} - u_{i-1/2}) \times u_{i-1/2} = - u_{i-1/2}\times u_{i+1/2} = -\sin(2\theta_i),
$$
as claimed.
\proofend

Next, one computes the commutators $[v_i,v_j]$. It is clear that
if $|i-j|\geq 2$ then $v_i$ and $v_j$ commute. One needs to
compute $[v_{j-1/2},v_{j+1/2}]$ for a half-integer $j$.

Let $w_j$ be the following tangent vector to ${\cal C}_n$: it
moves only vertices $O_{j-1/2}$ and $O_{j+1/2}$ with the
velocities
\begin{equation} \label{velo}
\frac{u_j}{\sin^2 \theta_{j-1/2}}\ \ {\rm and}\ \ \frac{u_j}{\sin^2 \theta_{j+1/2}},
\end{equation}
respectively, and changes only the radii $r_{j-1/2}$ and
$r_{j+1/2}$ with the rates
$$
\frac{\cos 2 \theta_{j-1/2}}{\sin^2 \theta_{j-1/2}}\ \ {\rm and}\ \ \frac{\cos 2 \theta_{j+1/2}}{\sin^2 \theta_{j+1/2}},
$$
respectively.

\begin{lemma} \label{comm}
The commutator $[v_{j-1/2},v_{j+1/2}]$ is proportional to $w_j$.
\end{lemma}

\proof The vector $[v_{j-1/2},v_{j+1/2}]$ moves vertices
$O_{j-1/2}$ and $O_{j+1/2}$ and changes the radii $r_{j-1/2}$ and
$r_{j+1/2}$, leaving everything else intact. The computation of
the action on the vertices is the same as for the Birkhoff
distribution done in \cite{BZ}. We use the result: it is an
infinitesimal motion of these two vertices along the line
$O_{j-1/2} O_{j+1/2}$ such that the signed perimeter length
remains zero.

\begin{figure}[hbtp]
\centering
\includegraphics[height=0.6in]{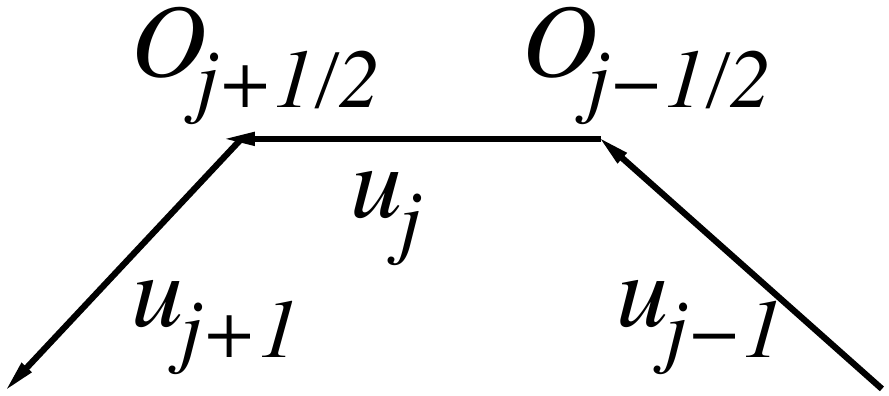}
\caption{To Lemma \ref{comm}.}
\label{speeds}
\end{figure}

Let us check that the velocities (\ref{velo}) satisfy this
description. Once again, different combinations of orientations
give the same result; let us consider the situation in Figure
\ref{speeds}. The rates of change of the three edge lengths are,
respectively,
\begin{equation} \label{threespeeds}
\frac{u_{j-1}\cdot u_j}{\sin^2 \theta_{j-1/2}} = \frac{\cos 2\theta_{j-1/2}}{\sin^2 \theta_{j-1/2}}, \frac{1}{\sin^2 \theta_{j+1/2}}- \frac{1}{\sin^2 \theta_{j-1/2}}, - \frac{u_{j+1}\cdot u_j}{\sin^2 \theta_{j+1/2}} = -\frac{\cos 2\theta_{j+1/2}}{\sin^2 \theta_{j+1/2}},
\end{equation}
and the sum is zero, as needed.

To calculate the rate of change of the radii, we use Lemma
\ref{signed} again. One has $|O_{j-1/2} O_{j-3/2}|=
r_{j-1/2}-r_{j-3/2}$, hence, using (\ref{threespeeds}),
$$
r_{j-1/2}'= \frac{\cos 2 \theta_{j-1/2}}{\sin^2 \theta_{j-1/2}}.
$$
Likewise, $|O_{j+1/2} O_{j+3/2}|= r_{j+3/2}-r_{j+1/2}$, hence
$$
r_{j+1/2}'= \frac{\cos 2 \theta_{j-1/2}}{\sin^2 \theta_{j+1/2}},
$$
as claimed.
\proofend

It is proved in \cite{BZ} that the projections of the $2n$ vectors
$v_i$ and $w_j$ to ${\cal E}_n$ span the tangent space at every
point. We need to show that the span of these vectors also
contains the fiber of the projection  $\pi: {\cal C}_n \to {\cal
E}_n$. To this end, we compute the kernel of $d\pi$.

\begin{lemma} \label{kern}
The vector
$$
\xi= \sum w_j + \sum \frac{\cos\theta_i}{\sin^3\theta_i}\ v_i
$$
generates Ker $d\pi$.
\end{lemma}

\proof We need to show that $\xi$ does not move any vertex.
Indeed, the velocity of $O_i$ corresponding to the tangent vector
$\xi$ is
\begin{equation} \label{vert}
\frac{\cos\theta_i}{\sin^3\theta_i} J(u_{i+1/2}-u_{j-1/2}) + \frac{1}{\sin^2\theta_i} (u_{i+1/2}+u_{j-1/2}).
\end{equation}
Without loss of generality, assume that
$$
u_{i+1/2} = (\cos\theta_i,\sin\theta_i),\ u_{i-1/2} = (\cos\theta_i,-\sin\theta_i),
$$
and it follows by a direct computation that (\ref{vert}) vanishes.

Finally, using the definitions of $w_j$ and $v_i$, we compute the
rate of change of the radii under $\xi$:
$$
r_i'=2\frac{\cos 2\theta_i}{\sin^2\theta_i} - \sin 2\theta_i \frac{\cos\theta_i}{\sin^3\theta_i} = -2.
$$
Thus $\xi$ corresponds to the deformation described in Lemma
\ref{radius}, adding a constant to all the radii. \proofend

This completes the proof of the theorem.
\proofend


\section{A distribution on ${\cal F}_n$ and $n$-fine curves} \label{distrF}

Recall the correspondence ${\cal C}_n^\circ \to {\cal F}_n^\circ$
from Proposition \ref{isom}. Given an oriented chain of circles
with centers $O_1,O_2,\ldots$ and signed radii $r_1, r_2,\ldots$,
let us now label the tangency points of the consecutive circles by
half-integers: $B_j$ is the tangency point of the circles with
centers $O_{j-1/2}$ and $O_{j+1/2}$ where $j\in \Z+\frac{1}{2}$. Then consistent orientations
of the circles define  directions $u_j$ at points $B_j$,
providing a framed polygon $(B_j,u_j)$.

As a result, we obtain an $n$-dimensional distribution ${\cal D}$ on ${\cal
F}_n^\circ$, the image of the distribution ${\cal \widetilde B}$ on ${\cal C}_n^\circ$.

Let us give a geometric description of the distribution ${\cal
D}$. Consider a framed $n$-gon $(B_j, u_j)$. If $n$ is odd,
consider the motion of vertex $B_k$ with the velocity $u_k$; all
other vertices do not move. By Lemma \ref{oddeven}, the framing
direction are determined by the polygon, and we obtain a tangent
vector to ${\cal F}_n$, denoted by $\eta_k$.

If $n$ is even, this approach does not work because moving only
one vertex of a polygon may violate the condition of Lemma
\ref{oddeven}. Instead, we describe the distribution as spanned by
certain motions of the sides of a polygon. Namely, for an integer
$m$, consider the  motion of two vertices, $B_{m-1/2}$ and
$B_{m+1/2}$, such that the velocity of $B_{m-1/2}$ is proportional
to $u_{m-1/2}$, the velocity of $B_{m+1/2}$ is proportional to
$u_{m+1/2}$, and the condition of  Lemma \ref{oddeven} holds. The
rest of the vertices and the framing vectors remain intact. This
motion defines a tangent vector to ${\cal F}_n$ which we denote by
$\nu_m$.

\begin{lemma} \label{corresp}
For odd $n$, the distribution ${\cal D}$ is generated by the
vectors $\eta_k$, and for even $n$, by the vectors $\nu_m$.
\end{lemma}

\begin{figure}[hbtp]
\centering
\includegraphics[height=2in]{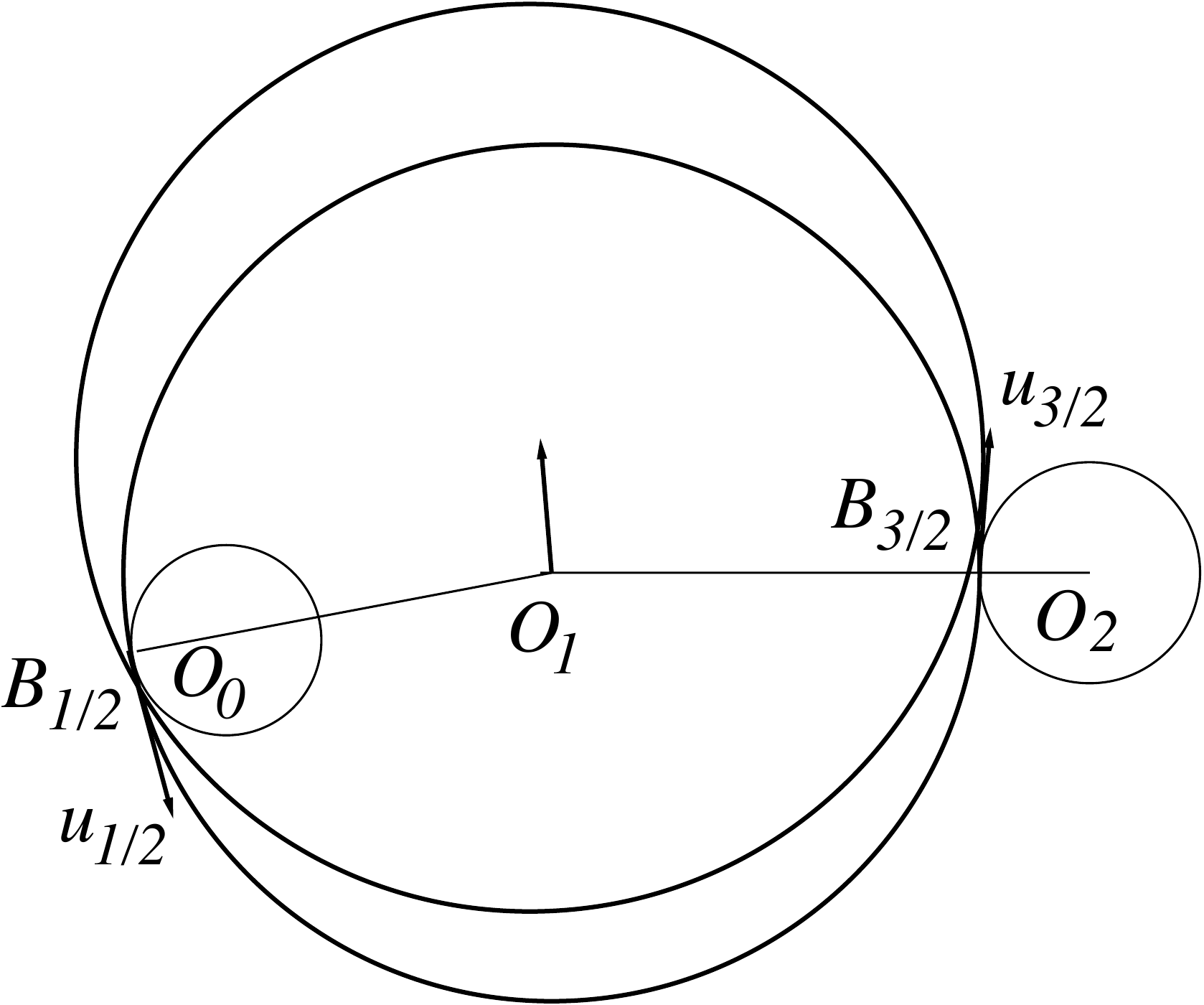}
\caption{Proof of Lemma \ref{corresp}.}
\label{deform}
\end{figure}

\proof Consider  a generator of the distribution ${\cal \widetilde
B}$ corresponding to an infinitesimal motion of the center $O_1$
along a hyperbola with foci $O_0$ and $O_2$ in Figure
\ref{deform}. The perturbed  circle remains tangent to the fixed
circles with centers $O_0$ and $O_2$, therefore the tangency
points $B_{1/2}$ and $B_{3/2}$ move along these circles, that is,
their velocities are proportional to $u_{1/2}$ and $u_{3/2}$. This
is precisely the description of the vector $\nu_1$.

For odd $n$, this tangent vector to ${\cal F}_n$ is a linear
combination of $\eta_{1/2}$ and $\eta_{3/2}$. Hence ${\cal D}$ is
a subspace of the space spanned by the vectors $\eta_k$. Since the
dimensions  are equal, these spaces coincide. \proofend

The distribution ${\cal D}$ comes equipped with a {\it positive
cone} consisting of those infinitesimal motions of framed polygons
for which each vertex $B_j$ moves with the velocity that is a
positive multiple of the  vector $u_j$.

The distribution ${\cal D}$ is defined on non-generic framed
polygons as well, but its bracket generating properties do not
extend to the whole of ${\cal F}_n$.

\begin{example} \label{circle}
{\rm Fix an oriented circle and consider the set of $n$-gons
inscribed in this circle and framed by the tangent vectors to it.
One obtains an $n$-dimensional submanifold $S_n\subset {\cal
F}_n$, tangent to ${\cal D}$, that is, $S_n$ is an $n$-dimensional
leaf of ${\cal D}$. Note that the respective chains of circles are
degenerate: all the circles in the chains coincide. }
\end{example}

Now we define a geometric object of the main interest of this study.

\begin{definition} \label{fine}
An $n$-fine curve is an oriented curve $\gamma$ with a 1-parameter
family of inscribed framed $n$-gons $B_1(t),\ldots,B_n(t)$, whose
framing is given by the positive unit tangent vectors of $\gamma$
at the respective points $B_j(t)$; here $t\in [0,1]$ is the
parameter in the family. The family satisfies the assumptions:
\begin{enumerate}
\item the velocity $B_j'(t)$ has the positive direction for all
values of $t$ and $j$ (that is, all vertices move along $\g$ with
positive speeds); \item the vertices of the polygon are cyclically
permuted over a period:
$$
B_j(t+1)=B_{j+1}(t)\  {\rm for\ all} \ t\  {\rm and}\   j;
$$
\item for $n\ge 3$, the tangent lines to $\g$ at points $B_j(t)$
and $B_{j+1}(t)$ are not parallel for all $t$ and $j$.
\end{enumerate}
\end{definition}

The 1-parameter families of inscribed polygons $B(t)$ are
considered up to reparameterizations. 

\begin{figure}[hbtp]
\centering
\includegraphics[height=1.5in]{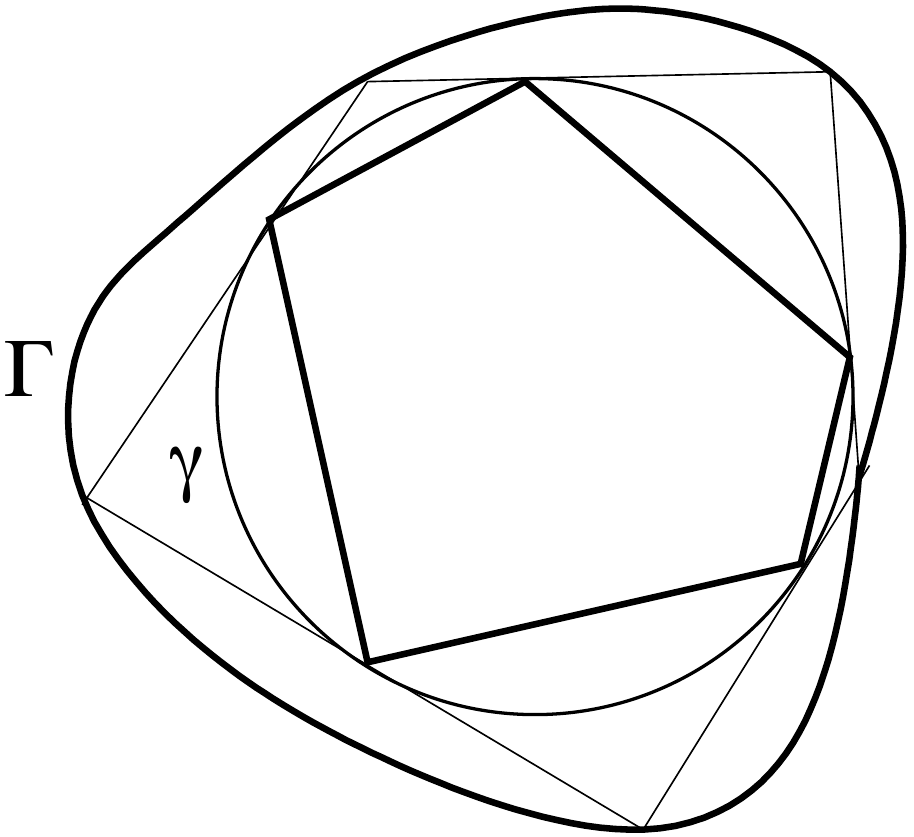}\qquad
\includegraphics[height=1.5in]{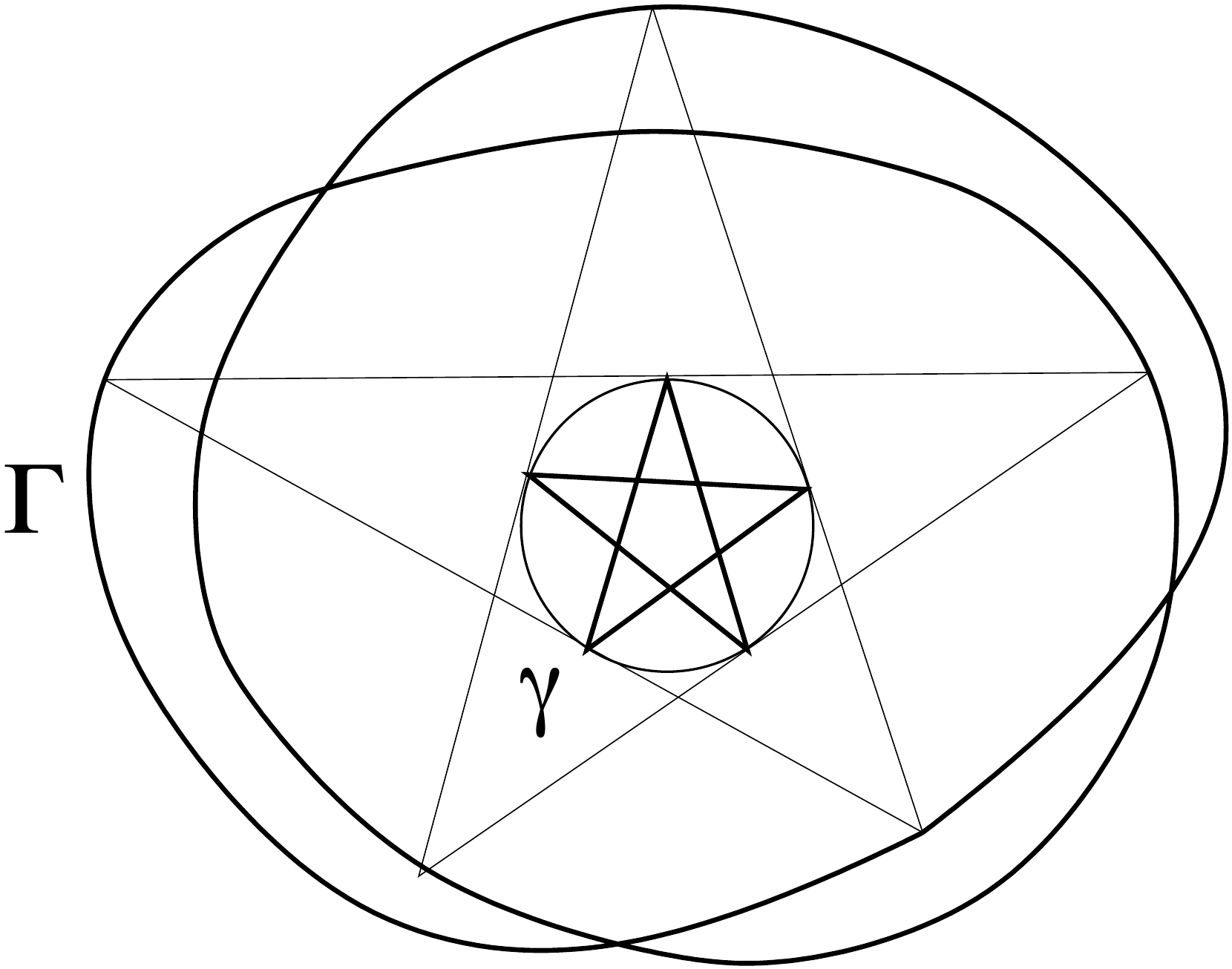}
\caption{Simple and star-shaped polygons in a $5$-fine curve.}
\label{stars}
\end{figure}

The locus of the intersection points of the tangent lines 
at points $B_j(t)$ and $B_{j+1}(t)$ is a curve $\G$ contained in the equitangent locus of the curve $\g$, see Figure \ref{stars}.

\begin{figure}[hbtp]
\centering
\includegraphics[height=1.2in]{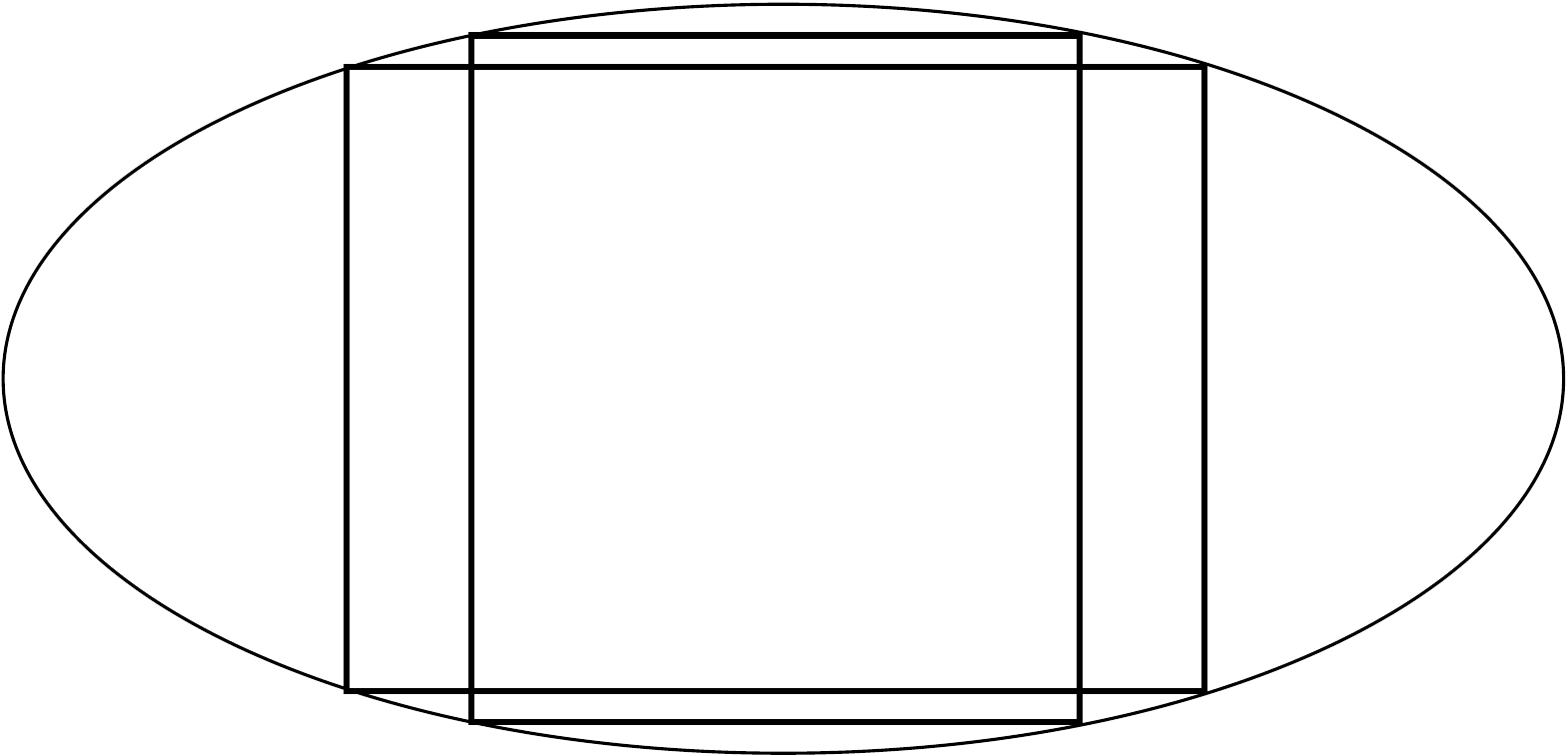}
\caption{An ellipse is not a $4$-fine curve.}
\label{ellipse2}
\end{figure}

\begin{example} \label{cirandell}
{\rm A circle provides an example of an $n$-fine curve for all
$n$. In fact, a circle gives infinitely many $n$-fine curves,
since any inscribed polygon is framed by the tangent directions to
the circle.

An ellipse is not a $4$-fine curve, in spite of the fact that at
admits a 1-parameter family of inscribed framed quadrilaterals:
not all the  vertices move in the ``right" direction, see Figure
\ref{ellipse2}. }
\end{example}

 A curve in a manifold
with a distribution is called {\it horizontal} if it is everywhere
tangent to the distribution. If a distribution is equipped with a
field of cones, a horizontal curve is called {\it positive} if its
velocity at every point belongs to the respective cone.

The cyclic group $\Z_n$ acts on ${\cal F}_n$ by cyclically
permuting the vertices of a polygon; this action preserves the
distribution ${\cal D}$. Denote by $\sigma$ the generator of this
group $(1,2,\ldots,n) \mapsto (2,3,\ldots,n,1).$

The above discussion implies the following result relating $n$-fine curves with
horizontal curves in $({\cal F}_n,{\cal D})$.

\begin{proposition} \label{horiz}
An $n$-fine curve $\g$ lifts to a positive horizontal curve
$\widetilde\g(t)=(B_1(t),\ldots,B_n(t))$ in ${\cal F}_n$
satisfying the monodromy condition
$\widetilde\g(t+1)=\sigma(\widetilde\g(t))$. Conversely, given
such a horizontal curve, the plane curves $B_j(t),\ t\in[0,1],
j=1,\ldots,n$, combine to form an $n$-fine curve.
\end{proposition}


\section{A case study: bigons} \label{chord}

In this section we consider the special case of framed bigons.

By a framed 2-gon we mean a segment $B_1 B_2$ with unit vectors
$u_1$ and $u_2$ attached to points $B_1$ and $B_2$, such that
$$
\angle (u_1, B_1 B_2) = \angle (B_1 B_2, u_2).
$$
Thus if $(B_1 B_2, u_1,u_2)$ is a framed 2-gon then so is $(B_2
B_1, u_2, u_1)$ (this is the action of $\Z_2$ on  ${\cal F}_2$).

\begin{figure}[hbtp]
\centering
\includegraphics[width=1.7in]{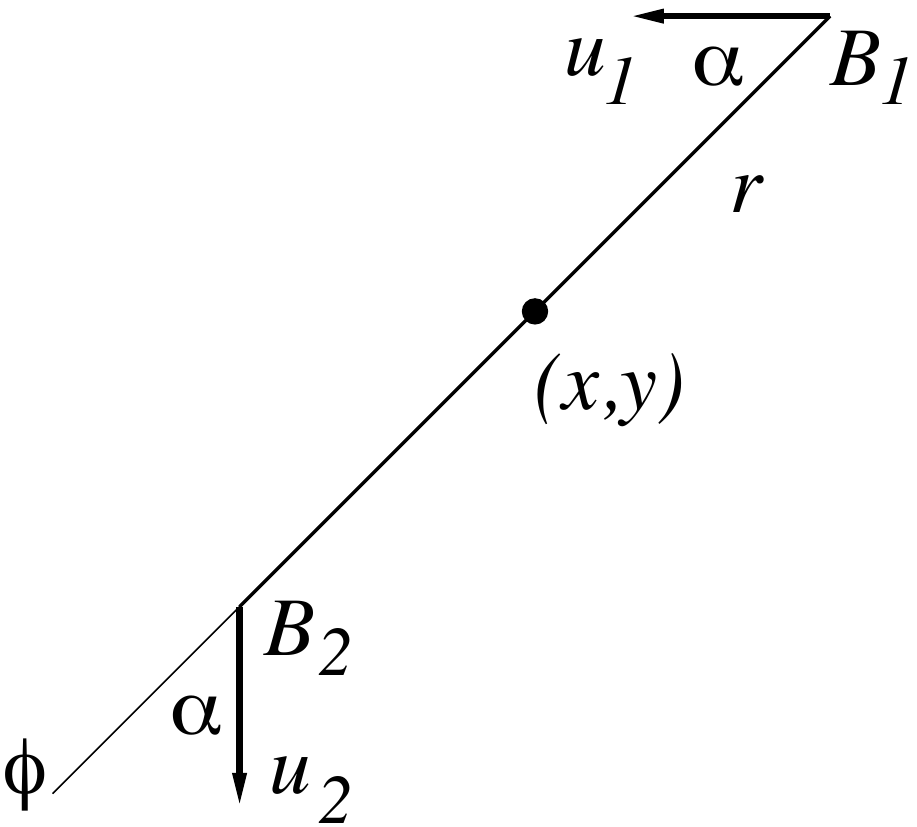}
\caption{A framed 2-gon.}
\label{segment}
\end{figure}

One can think of a framed 2-gon as a mechanical device consisting
of a telescopic axle of variable length with two wheels making the
same angle with the axle, see Figure \ref{segment}.

The space ${\cal F}_2$ is 5-dimensional: the general formula dim
${\cal F}_n = 2n$ does not apply since the condition of Lemma
\ref{oddeven} automatically holds for $n=2$. We introduce
coordinates as follows: $(x,y)$ is the midpoint of the segment
$B_1 B_2$, $\varphi$ is its direction, $r$ is its half-length, and
$\alpha$ is the angle made by the framing vectors with the
segment.

The distribution ${\cal D}$ is defined by the condition that the
velocities of $B_1$ and $B_2$ are aligned with $u_1$ and $u_2$,
respectively. The positive cone is defined by the condition that
these velocities are proportional to $u_1$ and $u_2$ with positive
coefficients.

Similarly to Theorem \ref{gen}, the distribution ${\cal D}$ is non-integrable.

\begin{theorem} \label{gen2}
The distribution ${\cal D}$ is bracket generating of the type $(3,5)$.
\end{theorem}

\proof
One has:
$$
B_1=(x-r\cos\varphi, y-r\sin\varphi),\ B_2=(x+r\cos\varphi, y+r\sin\varphi).
$$
The velocities of $B_1$ and $B_2$ are proportional to $u_1$ and
$u_2$, respectively. Differentiate $B_1$ and $B_2$, and take cross
product with
$$
u_1=(\cos (\varphi-\alpha), \sin (\varphi-\alpha)),\ u_2= (\cos (\varphi+\alpha), \sin (\varphi+\alpha)),
$$
to obtain, after some  transformations, the following two
relations:
\begin{equation*} \label{2rel}
x' \sin \varphi - y' \cos \varphi + r' \tan \alpha =0, \ x' \cos \varphi + y' \sin \varphi - r \varphi' \cot \alpha=0.
\end{equation*}
These relations define the distribution ${\cal D}$ in coordinates
$(x,y,r,\alpha,\varphi)$.

Change coordinates $(x,y)$ to $(p,q)$ as follows (this coordinate
change is suggested by the theory of support function in convex
geometry):
\begin{equation*}
x=p \sin\varphi - q \cos\varphi,\ y=-p \cos\varphi - q \sin\varphi.
\end{equation*}
In the new coordinates, the distribution ${\cal D}$ is given by
the equations
\begin{equation*}
p'+q\varphi'+r'\tan \alpha =0, \ p \varphi' - q' -  r \varphi' \cot \alpha=0.
\end{equation*}
It follows that ${\cal D}$ is defined by the 1-forms
\begin{equation} \label{forms}
\theta_1=dp+q\ d\varphi+\tan \alpha\ dr,\ \theta_2 = dq- (p-  r \cot \alpha)\ d\varphi,
\end{equation}
and is generated by the vector fields
$$
\nu=\frac{\partial}{\partial\alpha},\ \xi=\frac{\partial}{\partial\varphi}-q\frac{\partial}{\partial p}+(p-r \cot \alpha) \frac{\partial}{\partial q},\ \eta=\tan \alpha \frac{\partial}{\partial p} - \frac{\partial}{\partial r}.
$$
Therefore
$$
[\nu,\xi] = \frac{r}{\sin^2\alpha}\ \frac{\partial}{\partial q},\ [\nu,\eta] = \frac{1}{\cos^2\alpha}\ \frac{\partial}{\partial p}.
$$
This implies that the first order commutators of vector fields
tangent to ${\cal D}$ generate the tangent bundle of  ${\cal F}_2$
at every point. \proofend

\begin{remark}
{\rm A bracket generating distributions of the type $(3,5)$ contains a 2-dimensional sub-distribution that bracket generates it, see \cite{Mon}, sect. 6.8. Thus the study of distributions of types $(3,5)$ and $(2,3,5)$ is essentially the same. 
We do not know a geometrical interpretation of this 2-dimensional distribution in our setting. Let us also mention that distributions of the type $(3,5)$ arise in rolling a surface on another surface, for example, a ball rolling on a plane. 
}
\end{remark}

For $n=2$, we need to remove the last item in the definition of
$n$-fine curves. The issue is that if a segment $B_1 B_2$ makes
half a turn inside a curve $\g$, so that the endpoints $B_1$ and
$B_2$ swap positions, there must be a moment when the tangent
lines to $\g$ at the endpoints of the segment are parallel (this
readily follows from the intermediate value theorem).

Thus we do not exclude the case when  these tangents intersect at
infinity. Then the locus of the intersection points, $\Gamma$,
becomes a non-contractible curve in the projective plane. For
example, as was mentioned earlier, there is an
infinite-dimensional family of 2-fine curves $\g$ for which
$\Gamma$ is a straight line.

Now we shall show that 2-fine curves are very flexible and vary in
infinite-dimensional families. 
We shall use material from \cite{Mon},  especially chapter 5 and
appendix D.

Let $M$ be a manifold equipped with a $k$-dimensional distribution
${\cal D}$. One considers the space $H^2_{\cal D}$ of
parameterized horizontal curves $[0,1]\to M$ whose first two
derivatives are square integrable. For a point $x\in M$, the space
$H^2_{\cal D}(x)$ of curves starting at $x$ is a Hilbert manifold.
The functional dimension of this manifold equals $k$.

Now we consider $M={\cal F}_2$ with its distribution ${\cal D}$.

\begin{theorem} \label{flexbig}
Given a strictly convex $2$-fine curve $\g$, there exists a set of
$2$-fine curves sufficiently close to $\g$, which is in one-to-one
correspondence with a codimension 5 submanifold of the Hilbert
manifold $H^2_{\cal D}(x)$.
\end{theorem}

\proof The 2-fine curve $\g$ gives rise to a horizontal curve
$\widetilde\g(t)$ in ${\cal F}_2$, see Proposition \ref{horiz}. Let
$x=\widetilde\g(0), y =\widetilde\g(1)$; these endpoints
correspond to the same chord of $\g$ with the opposite
orientations. Consider the set of horizontal curves $H^2_{\cal
D}(x,y)\subset H^2_{\cal D}(x)$ with the terminal point $y$. A
horizontal curve from $H^2_{\cal D}(x,y)$, which is sufficiently
close to $\widetilde\g$, defines a 2-fine curve close to $\g$.

To describe a neighborhood of $\widetilde\g$ in $H^2_{\cal
D}(x,y)$, let $\pi: H^2_{\cal D}(x) \to {\cal F}_2$ be the map
that assigns the terminal point to a path. This map is smooth, and
if it is a submersion at $\widetilde\g$, then $H^2_{\cal
D}(x,y)=\pi^{-1} (y)$ is a submanifold in $H^2_{\cal D}(x)$ of
codimension equal to the dimension of ${\cal F}_2$, that is, of
codimension 5. A curve which is a singular point of the map $\pi$
is called singular.

To prove that $\widetilde\g$ is not singular, we use the following
criterion from \cite{Mon}. Given a manifold $M$ with a
distribution ${\cal D}$, let $\theta_i$ be a basis of 1-forms that
define ${\cal D}$. Consider the differential 2-forms
$\Omega(\lambda)=\sum \lambda_i d\theta_i$ (where $\lambda_i$ are
Lagrange multipliers, not all equal to zero). Then a horizontal
curve $\widetilde\g$ is singular if and only if its tangent vector
at every point lies in the kernel of $\Omega(\lambda)$ for some
choice of $\lambda$.

Now, to a calculation. According to (\ref{forms}),
$$
d\theta_1= dq\wedge d\varphi\ - \sec^2\alpha\ d\alpha\wedge dr,\ d\theta_2= d\varphi\wedge d(p+  r \cot \alpha).
$$
One has
$$
\Omega(\lambda) = \lambda_1 (dq\wedge d\varphi - \sec^2\alpha\ d\alpha\wedge dr) +
\lambda_2 (d\varphi\wedge dp + \cot\alpha\ d\varphi\wedge dr - r \csc^2\alpha\ d\varphi\wedge d\alpha).
$$
Let $\widetilde\g=(p(t),q(t),r(t),\varphi(t),\alpha(t))$, and set
$\widetilde\g' = w$. Assume that $i_w \Omega(\lambda)=0$. Since
$$
 i_w \Omega(\lambda) = \lambda_2 \varphi' dp - \lambda_1 \varphi' dq +\dots
 $$
 where dots mean a combination of $dr, d\alpha$ and $d\varphi$, and since
 both $\lambda_1$ and $\lambda_2$ cannot vanish simultaneosly, we conclude that $\varphi'=0$.

 However, if a $2$-fine curve is strictly convex, then the first item of
 the definition of $2$-fine curves implies that the chord $B_1(t) B_2(t)$
 rotates with a nonzero rate, that is, $\varphi'\neq 0$, a contradiction.
 \proofend

\begin{remark}
{\rm A classic example of a singular horizontal curve is the curve
$\gamma(t)=(t,0,0),\ t\in[0,1]$ tangent to the distribution
$dz=y^2dx$ in 3-space. If   $(x(t),y(t),z(t))$ is its perturbation
as a horizontal curve with fixed end-points then
$$
0=z(1)-z(0)=\int_0^1 y^2(t) x'(t) dt.
$$
Since $x'(t)>0$, we have $y(t)=0$, and hence $z(t)=0$, for all
$t$; thus the perturbed curve is a reparameterization of $\gamma$.
Therefore the curve $\g$ is rigid, see \cite{Mon}. }
\end{remark}


\section{Triangles and quadrilaterals} \label{tri}

Let us consider $n$-fine curves for small values of $n$, namely,  $n=3$ and $n=4$.

\begin{proposition} \label{n=3}
A 3-fine curve is a circle.
\end{proposition}

\proof Any triangle is inscribed into a circle and, as a  framed
polygon, is of the kind described in Example \ref{circle}. Hence
the corresponding horizontal curve in ${\cal F}_3$ stays in $S_3$,
that is, the respective triangles are inscribed into the same
circle. \proofend

\begin{remark}
{\rm The distribution ${\cal D}$ on ${\cal F}_3$ is integrable,
that is, it is a foliation. This foliation consists of the level
surfaces of three functions: the radius and the two components of
the center of the incircle. This does not contradict Theorem \ref{gen} since ${\cal F}_3^{\circ}=\emptyset$.}
\end{remark}

A similar rigidity results holds for  $n=4$.

\begin{proposition} \label{n=4}
A  $4$-fine curve is a circle.
\end{proposition}

\proof As we mentioned earlier, a framed quadrilateral is cyclic. 

\begin{figure}[hbtp]
\centering
\includegraphics[height=1.7in]{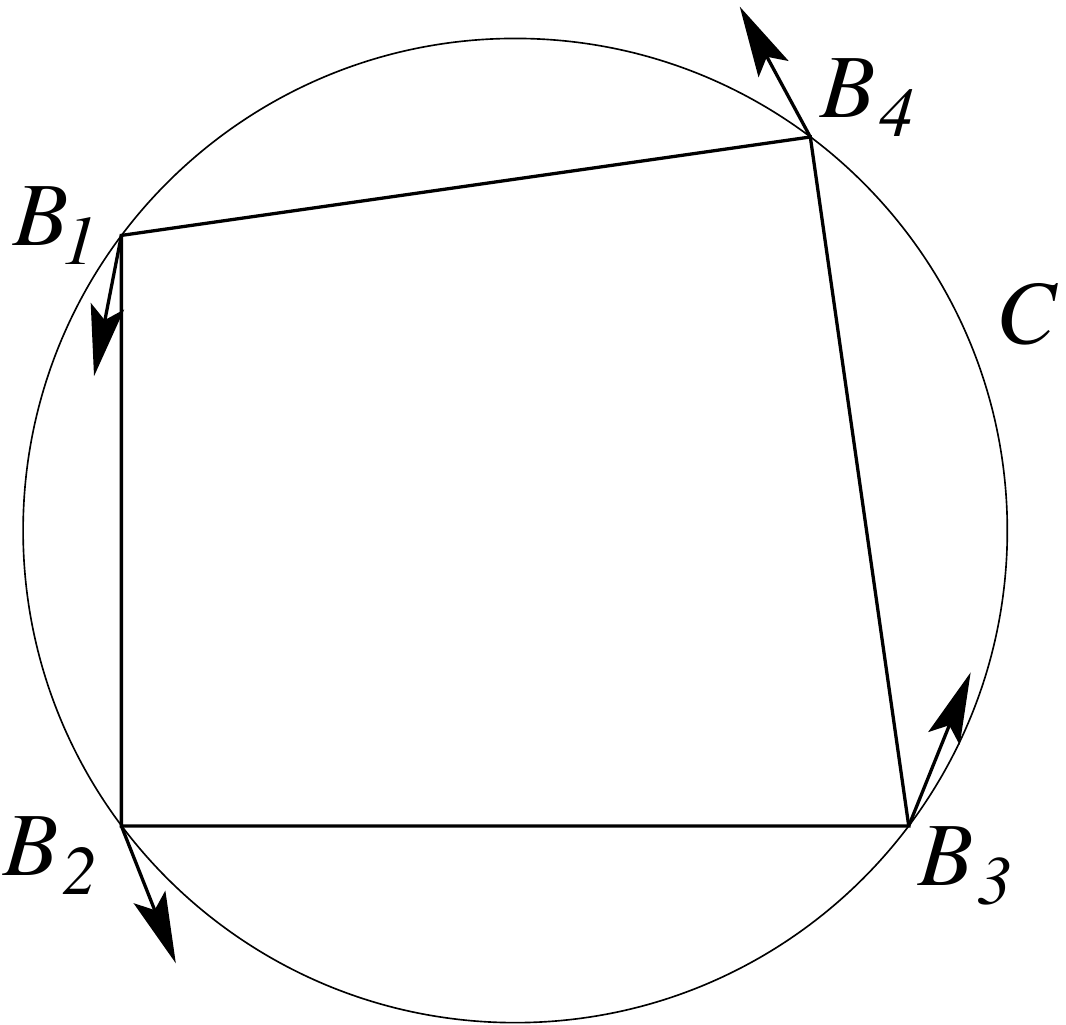}
\caption{Proof of Proposition \ref{n=4}}
\label{quad}
\end{figure}

According to Lemma \ref{oddeven}, the framing of our inscribed
quadrilateral $P=(B_1B_2B_3B_4)$ is obtained from the tangent
vectors to the circumscribed circle, say, $C$, by rotating the odd
and the even ones the same amount, say, $\varphi$, in the opposite
directions. See Figure \ref{quad}.

Let $\gamma$ be a  $4$-fine curve into which $P$ is inscribed.

If $\varphi=0$ then we are in the situation encountered earlier,
for $n=3$. Suppose that $\varphi > 0$. Then $\gamma$ intersects
the circle transversally.

Let $P'$ be a quadrilateral,  infinitesimally close to $P$,
inscribed into $\gamma$. Then $P$ is also inscribed into a circle,
$C'$. Note that the odd vertices of $P'$ are inside $C$ and the
even ones are outside of it. Therefore the circles $C$ and $C'$
intersect in four points, a contradiction. \proofend

\begin{remark}
{\rm The situation with small values of $n$ --  flexibility for
$n=2$ and rigidity for $n=3,4$ -- resembles another problem of a
somewhat similar flavor, the 2-dimensional Ulam's problem. The
Ulam problem is to describe homogeneous bodies that float in
equilibrium in all positions \cite{Mau}. In dimension three and
higher, almost nothing is known, but in dimension two, a variety
of results is available (although a complete solution is missing).
A relevant quantity here is the perimetral density of the body,
i.e., if the waterlines divide the boundary of a body $K$ in a
constant ratio $\sigma:(1-\sigma)$  the $\sigma$
is called the \emph{perimetral density}. If this density is $1/2$, there
is a functional family of solutions \cite{Au}, but if the density
is $1/3$ or $1/4$, then the only solutions are round discs. See
\cite{BMO1,BMO2,Odan,Ta1,Va,We1,We2} for references.}
\end{remark}


\section{An example} \label{exa}

Here we construct a nested pair of curves $\g \subset \G$ such that $\G$ is contained in the equitangent locus of $\g$. Our approach is similar to that of \cite{Ta2}.

We start with a discrete version of the problem, where a smooth
curve is replaced by a polygon. For a convex polygon, the tangent
direction at an interior point of a side is the direction of this
side, and a tangent direction  at a vertex is that of a support
line at this vertex (thus there are many tangent directions at a
vertex). We describe a position of a framed 2-gon by a triple whose first
element is the chord, the second element is a tangent line at the
first end-point, and the third element a tangent line at the
second end-point.

\begin{figure}[hbtp]
\centering
\includegraphics[width=2.7in]{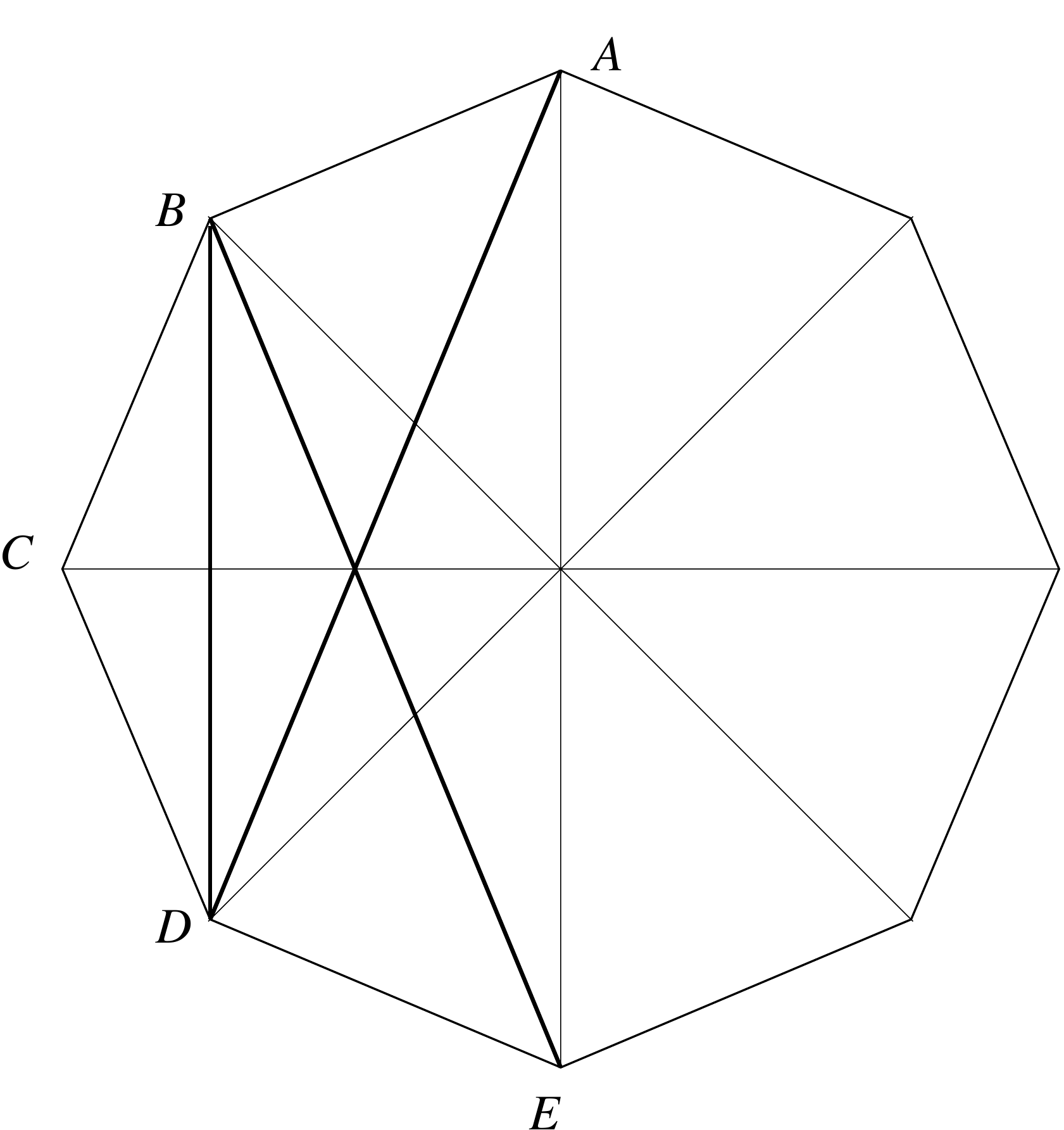}
\caption{Rotating a chord in a regular octagon}
\label{octa}
\end{figure}

Consider a regular octagon in figure \ref{octa}. Here is a
two-step continuous motion of the framed 2-gon $AD$ to $BE$:
$$
(AD, AB, DC) \mapsto (BD, AB, ED) \mapsto (BE, BC, ED).
$$
In the first step, as point $A$ moves toward $B$, the support
direction at the other end point $D$ turns accordingly. In the
second step, as point $D$ moves toward $E$, the support direction
at point $B$ turns accordingly. The whole process consists of 16
steps, so that the chord makes a complete circuit.

Next we approximate the octagon by a strictly convex smooth curve.
For example, we may replace each vertex by an arc of a circle of
very small radius and each side by an arc of a circle of very
large radius, thus approximating the octagon by a piecewise
circular curve. In this way one obtains a $C^1$-smooth 2-fine
curve $\g$ having dihedral symmetry. The outer curve $\Gamma$ is
the locus of intersection points of the respective pairs of
tangent lines of $\g$.

\begin{remark}
{\rm If $\g$ is a piecewise circular curve then $\Gamma$ is still
a polygonal line. Indeed, the equitangent locus of  two circles is
a line, the radical axis of the two circles. }
\end{remark}

\begin{figure}[hbtp]
\centering
\includegraphics[width=2.7in]{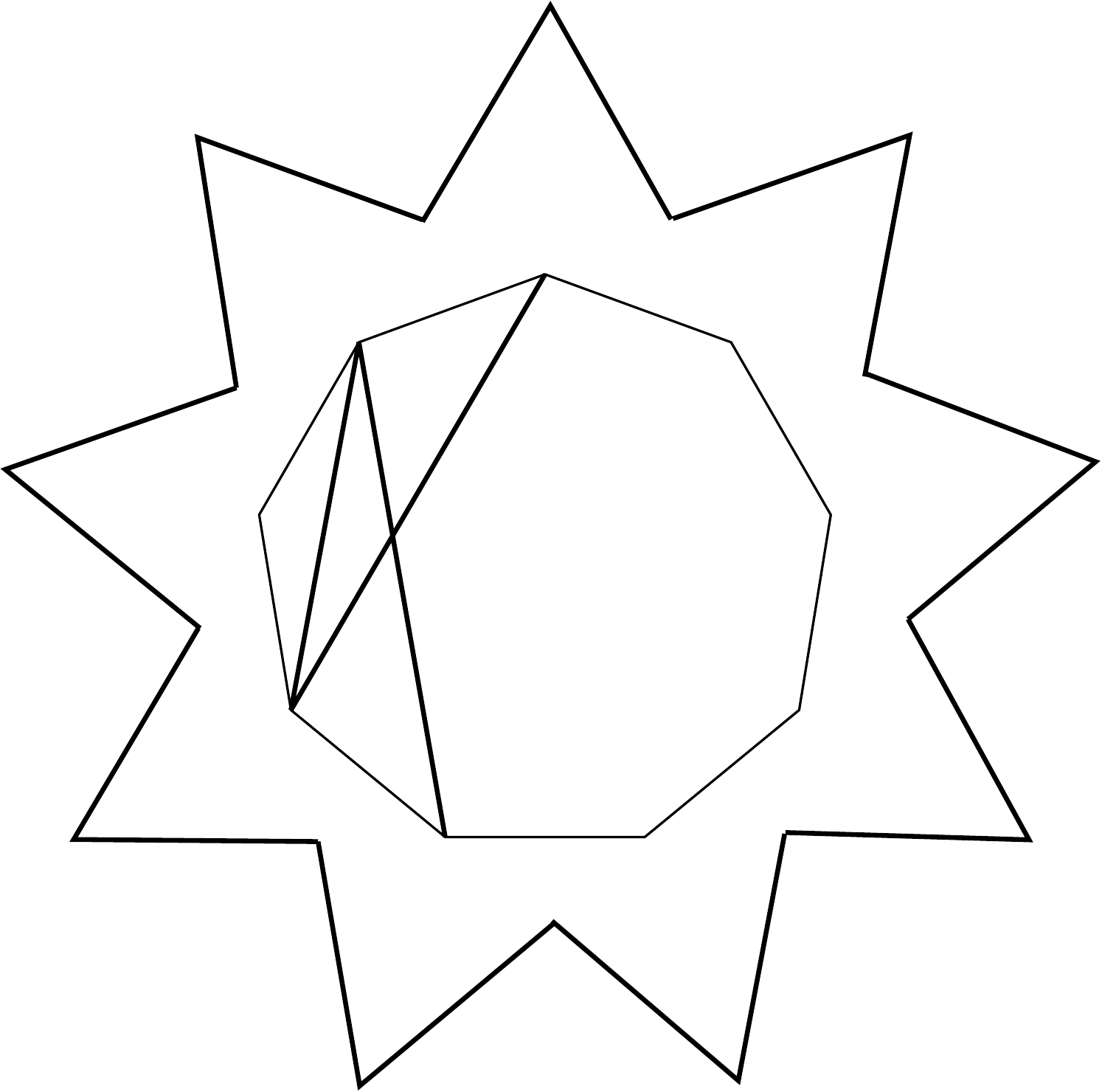}
\caption{Construction based on a regular nonagon.}
\label{nona}
\end{figure}

A similar construction can be made for other regular $n$-gons with
$n\ge 7$. See Figure \ref{nona} for $n=9$.

\section{Can the equitangent locus of an $n$-fine curve be a circle?} \label{canit}

As we mentioned earlier, the equitangent locus of a 2-fine curve may contain a straight line. 
This motivates the question in the title of this section. 

Taken literally, this question has an affirmative answer: one may
take $\g\subset\Gamma$ to be a pair of circles that support a family
of Poncelet $n$-gons, inscribed in $\Gamma$ and circumscribed
about $\g$; such polygons are called bicentric. 

The variety of
circles inside $\Gamma$ is 3-dimensional, and the existence of
bicentric $n$-gons describes a 2-dimensional subvariety. Denote
the radii of the circles by $r$ and $R$, and the distance between
their centers by $d$. Explicit formulas are available for the relation between these quantities (for $n=3$, due to
L. Euler, and for $n=4,5,6,7,8$, to N. Fuss), for example,
$$
R^2= d^2 + 2rR\ \ {\rm for}\ n=3, \ (R^2 -d^2)^2 =2r^2(R^2 +d^2)\ \ {\rm for}\ n=4,
$$
see, e.g., \cite{Ra} for a recent treatment.

Thus the real question is whether there exist a non-circular $n$-fine curve $\g$
whose equitangent locus contains a circle $\G$.

\begin{figure}[hbtp]
\centering
\includegraphics[width=2in]{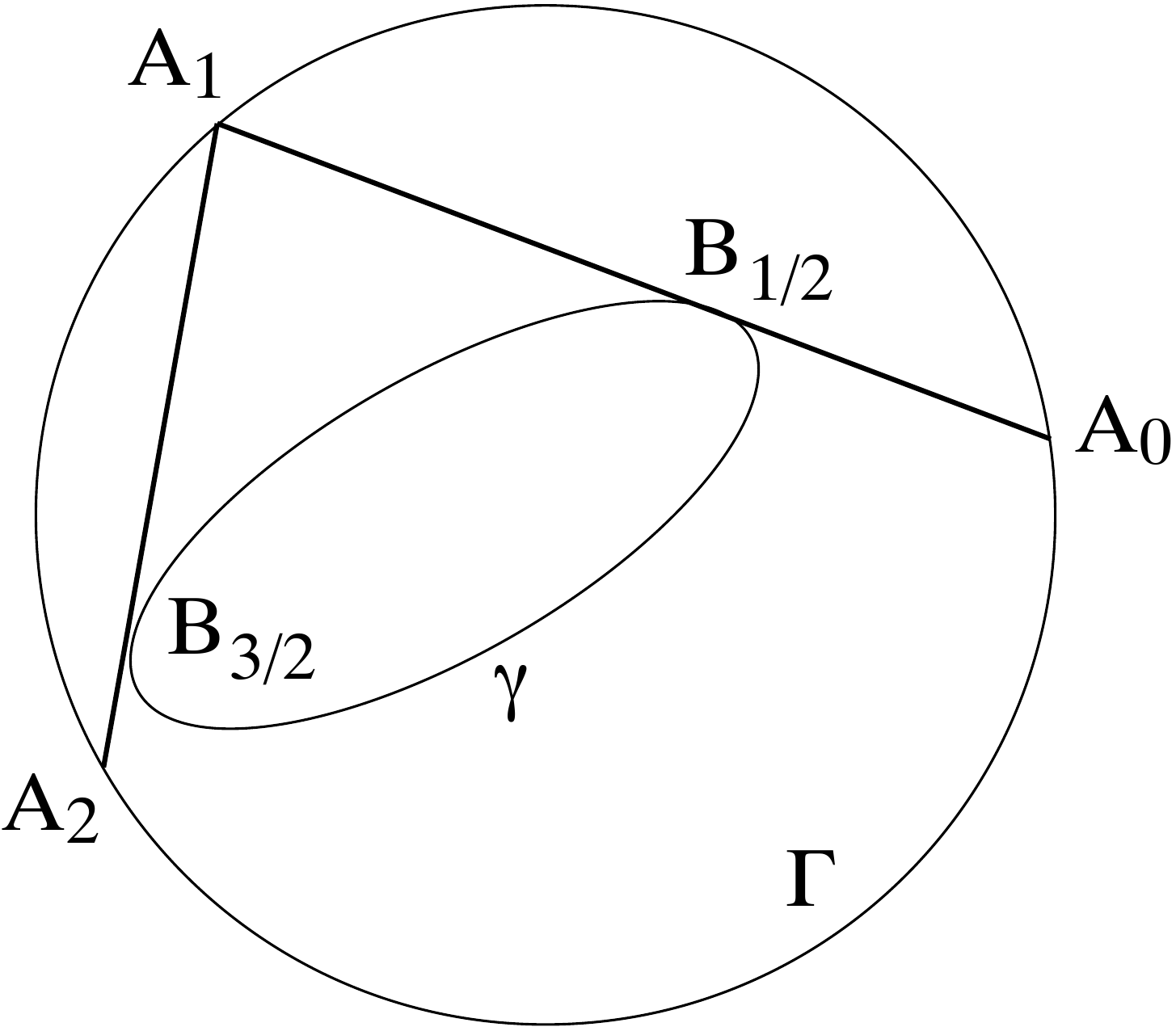}
\caption{An equitangent polygon inscribed into a circle.}
\label{setfield}
\end{figure}

Let $\Gamma$ be a unit circle; denote by ${\cal S}_n$ the set of
inscribed $n$-gons. Assume that $(A_1,A_2,\ldots)\in {\cal S}_n$
is an equitangent polygon, that is, there exists a  curve $\g$
such that the two tangent segments to $\g$ from each point $A_i$
are equal. Denote the tangency points by $B_j$, where $j$ is
half-integer, as shown in Figure \ref{setfield}.

The following lemma is an analog of Lemma \ref{oddeven}.

\begin{lemma} \label{oddeven1}
If $n$ is odd then the polygon $A$ uniquely determines the
positions of points $B_j$. If $n$ is even, then a necessary and
sufficient condition for $A$ to be an equitangent polygon is that
$$
\sum_{i=1}^n (-1)^i |A_i A_{i+1}| =0.
$$
If this relation holds then there exists a 1-parameter family of
the respective polygons $B$.
\end{lemma}

\proof Let $x_i=|A_{i} B_{i-1/2}|= |A_{i} B_{i+1/2}|$. Then one
has a system of linear equations
\begin{equation} \label{systl}
x_i + x_{i+1}=|A_i A_{i+1}|,\ i=1,\ldots,n.
\end{equation}
If $n$ is odd, this system has a unique solution, and if $n$ is
even, its matrix has rank $n-1$. The alternating sum of the
equations yields the relation claimed in the lemma.

We note that not every solution of the system is `geometric': for
example, if $x_i >  |A_i A_{i+1}|$ then the point $B_{i+1/2}$ lies
on an extension of the side $A_i A_{i+1}$. \proofend

Note that, for $n=4$, the condition of Lemma \ref{oddeven1} means that the respective quadrilateral is bicentric. 

To fix ideas,  assume that $n$ is odd: $n=2m+1$. We also assume
that the polygons $A$ are convex.

Consider an equitangent $n$-gon, as in Figure \ref{setfield}. Move
the vertex $A_1$ slightly. Then the vertex $A_2$ also moves, so
that the line $A_1A_2$ remains tangent to the curve $\g$, which
forces the vertex $A_3$ move as well, etc. The next lemma shows
that this process describes a well defined motion of an inscribed
$n$-gon.

\begin{lemma} \label{close}
In this process, the $(n+1)$st vertex moves the same amount as the
1st one, that is, the periodicity condition $A_{n+1}=A_1$ is
preserved.
\end{lemma}

\proof Let $A_i'$ be a perturbed vertex. By elementary geometry,
the infinitesimal triangles $A_i A_i' B_{i+1/2}$ and $A_{i+1}
A_{i+1}' B_{i+1/2}$ are similar, hence
\begin{equation} \label{Newteq}
\frac{|A_i A_i'|}{x_i}=\frac{|A_{i+1} A_{i+1}'|}{x_{i+1}}.
\end{equation}
Continuing this chain of equalities, and using the fact that
$x_{n+1}=x_1$, we find that $|A_1 A_1'|=|A_{n+1} A_{n+1}'|$, as
claimed. \proofend

\begin{figure}[hbtp]
\centering
\includegraphics[width=2in]{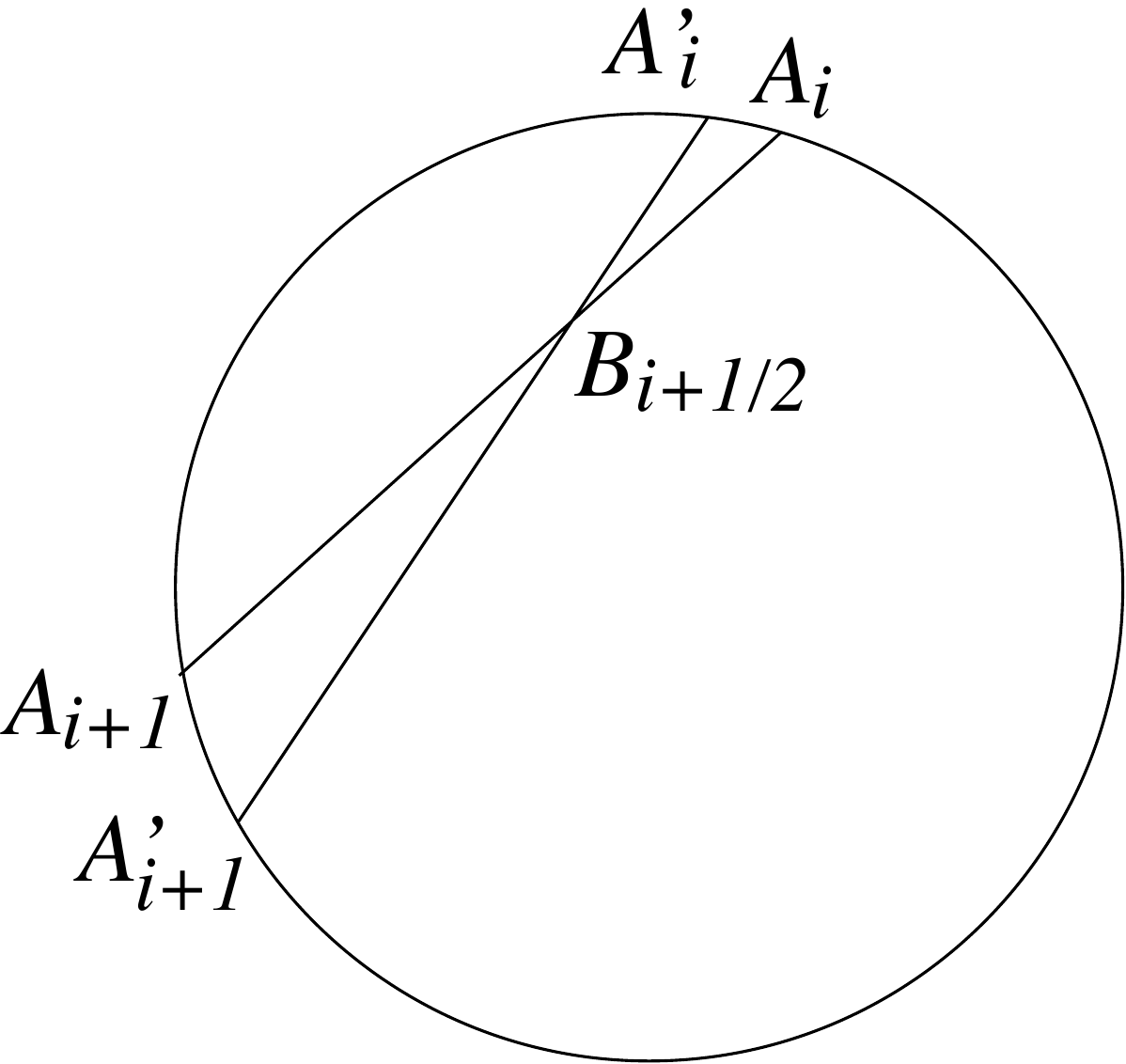}
\caption{Proof of Lemma \ref{close}.}
\label{Newt}
\end{figure}

Thus our construction provides a field of directions on the space of inscribed $n$-gons,  ${\cal
S}_n$. Assuming that vertex $A_1$ moves with unit speed, we obtain
a vector field on ${\cal S}_n$, which we denote by $\xi$.

Recall that $\sigma$ denotes the cyclic permutation of the
vertices of an $n$-gon. If a trajectory of $\xi$ connects a
polygon $A$ with the polygon $\sigma(A)$, then the envelopes of
the lines $A_i A_{i+1}$ combine to form the desired $n$-fine curve $\g$. For example, this is the case if $A$ is bicentric.

The case $n=3$ is trivial: every triangle is circumscribed about a
circle, hence the trajectories of $\xi$ on ${\cal S}_3$ consists
of the Poncelet families of triangles. The first non-trivial case
is that of pentagons.

We do not have a description of the trajectories of the vector
field $\xi$ on ${\cal S}_n$. As a first step toward this goal, we
consider a trivial solution, corresponding to a regular
$n$-gon (which is bicentric), and linearize the system near this solution.

Let $A_i=(\cos\alpha_i,\sin\alpha_i)$; we use $\alpha_i,\ i=1,\ldots,n$, as coordinates in ${\cal S}_n$.

\begin{lemma} \label{syst}
After a change of the time parameter, the vector field $\xi$
corresponds to the system of differential equations
$$
\dot\alpha_i = \sum_{k=0}^{n-1} (-1)^k \sin \frac{\alpha_{i+1+k}-\alpha_{i+k}}{2},\ i=1,\ldots,n.
$$
\end{lemma}

\proof
One has:
$$
|A_i A_{i+1}|=2\sin \frac{\alpha_{i+1}-\alpha_i}{2},
$$
and the solution to the system of  equations (\ref{systl}) is
$$
x_i = \sum_{k=0}^{n-1} (-1)^k \sin \frac{\alpha_{i+1+k}-\alpha_{i+k}}{2}.
$$
Using (\ref{Newteq}), we  change the time in our system so that
the vector field $\xi$ becomes $\dot\alpha_i= x_i$, and this
yields the result. \proofend

The solution corresponding to the regular $n$-gon is
\begin{equation*}
\alpha_i(t)=t \sin\frac{\pi}{n}  + \frac{2\pi i}{n},\ i=1,\ldots,n,
\end{equation*}
with  period
$$
T_0=\frac{2\pi}{\sin\frac{\pi}{n}}.
$$
Consider its infinitesimal perturbation:
\begin{equation} \label{perturb}
\alpha_i(t)=t \sin\frac{\pi}{n}  + \frac{2\pi i}{n} + \eps\beta_i(t).
\end{equation}

\begin{lemma} \label{linearized}
The linearized system of differential equations is
\begin{equation} \label{linsyst}
\dot\beta_i(t) = \cos\frac{\pi}{n} \left(\sum_{k=1}^{n-1} (-1)^{k-1} \beta_{i+k}\right),\ i=1,\ldots,n.
\end{equation}
\end{lemma}

\proof
Computing modulo $\eps^2$, we have:
$$
\sin \frac{\alpha_{i+1+k}-\alpha_{i+k}}{2} =  \sin\frac{\pi}{n} + \eps \cos\frac{\pi}{n} \left(\frac{\beta_{i+1+k}-\beta_{i+k}}{2}\right),
$$
hence
$$
x_i=\sin\frac{\pi}{n} + \eps \cos\frac{\pi}{n} \left(\sum_{k=1}^{n-1} (-1)^{k-1} \beta_{i+k}\right).
$$
Therefore the equation from Lemma \ref{syst} becomes (\ref{linsyst}).
\proofend

We note  an integral of the system (\ref{linsyst}): $\sum
\beta_i(t)=$ const. Since adding a constant to all $\beta_i$
amounts to a rotation of the circle, without loss of
generality, we assume that $\sum \beta_i(t)=0$.

The system has the circulant matrix
$$
M_n=\cos\frac{\pi}{n}\
\left(
\begin{array}{cccccc}
0&1&-1&1&\dots&-1\\
-1&0&1&-1&\dots&1\\
\dots&\dots&\dots&\dots&\dots&\\
1&-1&1&-1&\dots&0
\end{array}
\right).
$$
Let $\lambda_0=0,\lambda_1,\ldots,\lambda_{n-1}$ be its
eigenvalues.

\begin{lemma} \label{spec}
One has
$$
\lambda_j=-\lambda_{n-j}=\sqrt{-1} \cos\frac{\pi}{n} \tan\frac{\pi j}{n},\ \ j=1,\ldots,m.
$$
\end{lemma}

\proof The eigenvalues of a circulant matrix of order $n$ are the
values at $\omega^j$ of the polynomial whose coefficients are the
entries of the first row; here $\omega=\exp{(2\pi i/n)}$. In the
case at hand,
$$
\lambda_j=\cos\frac{\pi}{n}\ (\omega^j-\omega^{2j}+ \dots - \omega^{(n-1)j}) = \cos\frac{\pi}{n} \ \frac{\omega^{j/2}-\omega^{-j/2}}{\omega^{j/2}-\omega^{-j/2}} = i \cos\frac{\pi}{n} \ \tan\frac{\pi j}{n},
$$
as claimed.
\proofend

The phase space of the differential equation (\ref{linsyst})
splits into the sum of a one-dimensional subspace, corresponding
to the zero eigenvalue, and $m$ two-dimensional subspaces,
corresponding to the eigenvalues $\lambda_j$. The dynamics in the
1-dimensional subspace is trivial, and  in the 2-dimensional
subspaces, one has rotations with the angular velocities
$|\lambda_j|$.

\begin{example}
{\rm If $n=3$, all solutions are periodic, with the same period.
This corresponds to the fact that every triangle is included into
a 1-parameter family of Poncelet triangles.

If $n=5$, the ratio of the absolute values of the eigenvalues
equals $\sqrt{5}-2$. Therefore the eigenvalues are rationally
independent, and a generic trajectory is an irrational line on a
2-dimensional torus, see, e.g., \cite{Ar}: the only periodic
trajectories are the ones in the 2-dimensional eigenspaces, corresponding to convex and star-shaped regular pentagons. }
\end{example}

\begin{conjecture} \label{indep}
For odd $n\ge 7$, the eigenvalues $\lambda_j,\ j=1,\ldots,m$, are
rationally independent.
\end{conjecture}

Assuming this conjecture, we have an infinitesimal rigidity result.

\begin{proposition} \label{rigid} For odd $n$, the only infinitesimal deformations of the 1-parameter family of the regular $n$-gons in the class of equitangent $n$-gons inscribed into a unit circle, are the families of 
bicentric $n$-gons.
\end{proposition}

\proof Assume that the perturbed system (\ref{linsyst}) has a
periodic solution with period $T=T_0+\eps Q$. Then
$$
(t+T) \sin\frac{\pi}{n} + \eps \beta_i(t+T)=2\pi+t \sin\frac{\pi}{n} + \eps \beta_i(t)
$$
for all $t$ and $j$, and hence
\begin{equation} \label{diff}
\beta_i(t+T_0)-\beta_i(t)=-Q\sin\frac{\pi}{n}
\end{equation}
(as before, computing modulo $\eps^2$). Adding (\ref{diff}) for
$i=1,\ldots,n$, and using the fact that $\sum \beta_i=0$
identically, we conclude that $Q=0$. Hence the solution to
(\ref{linsyst}) is periodic with period
$T_0=2\pi/\sin\frac{\pi}{n}$.

Note that $|\lambda_1|=\sin\frac{\pi}{n}$, and the respective period
is $T_0$. If a  solution has non-zero components in at least two
eigenspaces, then  Conjecture \ref{indep} implies that the motion
is not periodic. Finally, if only one component is non-zero, and
it corresponds to $\lambda_j$ with $j>0$, then, due to  Conjecture
\ref{indep}, the period is not commensurable with $T_0$.

We conclude that the only infinitesimal deformations of the
regular $n$-gons belong to a 2-parameter family, corresponding to the eigenvalue 
$\lambda_1$ (the family is 2-parameter, and not 3-parameter,
because we factored out the rotations of the circle). Since the families of bicentric
$n$-gons are solutions and they also constitute a 2-parameter space of solutions, all
infinitesimal perturbations are bicentric. \proofend

We finish with another conjecture.

\begin{conjecture} \label{totalrig}
If $\g$ is an $n$-fine curve and the  vertices of the respective equitangent $n$-gons lie on a circle then $\g$ must be a circle as well.
\end{conjecture}


\begin{thebibliography}{99}

\bibitem{Ar} V. Arnold. {\it Ordinary differential equations.} Springer-Verlag, Berlin, 2006.

\bibitem{Au} H. Auerbach. {\it Sur un problem de M. Ulam concernant l'equilibre des corps flottants.} Studia Math. 7 (1938), 121--142.

\bibitem{BG} T. Banchoff, P. Giblin. {\it On the geometry of piecewise circular curves.} Amer. Math. Monthly 101 (1994),  403--416.

\bibitem{BZ} Yu. Baryshnikov, V. Zharnitsky. {\it Sub-Riemannian geometry and periodic orbits in classical billiards.} Math. Res. Lett. 13 (2006),  587--598.

\bibitem{BMO1} J. Bracho, L. Montejano, D. Oliveros. {\it A classification theorem for Zindler carrousels}. J. Dynam. and Control Syst. 7 (2001), 367--384.

\bibitem{BMO2} J. Bracho, L. Montejano, D. Oliveros. {\it Carousels, Zindler curves and the floating body problem.} Period. Math. Hungar. 49 (2004),  9--23.

\bibitem{BGG} J. Bruce, P. Giblin, C. Gibson. {\it Symmetry sets.} Proc. Roy. Soc. Edinburgh Sect. A 101 (1985),  163--186.

\bibitem{Da} Ph. Davis. {\it Circulant matrices.}   John Wiley \& Sons, 1979.

\bibitem{Fl} L. Flatto. {\it  Poncelet's theorem.} AMS, Providence, RI, 2009.

\bibitem{GB} P. Giblin, S. Brassett. {\it Local symmetry of plane curves.}  Amer. Math. Monthly 92 (1985),  689--707.

\bibitem{GT} D. Genin, S. Tabachnikov. {\it On configuration spaces of plane polygons, sub-Riemannian geometry and periodic orbits of outer billiards}. J. Mod. Dyn. 1 (2007),  155--173.

\bibitem{Gl1} A. Glutsyuk. {\it On odd-periodic orbits in complex planar billiards}.  
J. Dyn. Control Syst. 20 (2014), 293--306.

\bibitem{Gl2} A. Glutsyuk. {\it On 4-reflective complex analytic planar billiards.} Preprint arXiv:1405.5990.

\bibitem{JR} J. Jeronimo-Castro, E. Rold\'an-Pensado. {\it A characteristic property of the Euclidean disc}. Period. Math. Hungar. 59 (2009), 215--224.

\bibitem{JRT} J. Jeronimo-Castro, G. Ruiz-Hernandez, S. Tabachnikov.
{\it The equal tangents property}. Adv. Geom. 14 (2014), 447--453.

\bibitem{La} J. Landsberg.
{\it Exterior differential systems and billiards}. Geometry, integrability and quantization, 35--54, Softex, Sofia, 2006.

\bibitem{Mau} R. Mauldin. {\it The Scottish book, Mathematics from the Scottish Caf\'e}. Birhauser, 1981.

\bibitem{Mon} R. Montgomery.
\textsl{A tour of subriemannian geometries, their geodesics and applications},
Amer. Math. Soc., Providence, RI, 2002.

\bibitem{Odan} K. Odani. {\it On Ulam's floating body problem of two dimension}.
Bulletin of Aichi Univ. of Education 58 (Natural Sciences) (2009),
1--4.

\bibitem{RT} H. Rademacher,  O. Toeplitz. \emph{The enjoyment of mathematics.} Princeton University Press, 1957.

\bibitem{Ra} M. Radi\'c. {\it One system of equations concerning bicentric polygons.} J. Geom. 91 (2009),  119--139.

\bibitem{Ta1} S. Tabachnikov. {\it Tire track geometry: variations on a theme}. Israel J. Math., 151 (2006), 1--28

\bibitem{Ta2} S. Tabachnikov. {\it The (un)equal tangents problem}.  Amer. Math. Monthly 110 (2012), 398--405.

\bibitem{TZ} A. Tumanov, V.  Zharnitsky. {\it Periodic orbits in outer billiard. } Int. Math. Res. Not. 2006, Art. ID 67089, 17 pp.

\bibitem{Va} P. V\'arkonyi. {\it Floating body problems in two dimensions.} Stud. Appl. Math. 122 (2009),  195--218.

\bibitem{We1} F. Wegner. {\it Floating bodies of equilibrium}. Studies Appl. Math. 111 (2003), 167--183.

\bibitem{We2} F. Wegner. {\it Three problems, one solution}. A web site. \url{http://www.tphys.uni-heidelberg.de/~wegner/Fl2mvs/Movies.html}

\end{thebibliography}
\end{document}